\documentclass[12pt]{article}
\usepackage{amssymb}
\usepackage{amsmath}
\usepackage{amsthm}
\usepackage{tikz}
\usepackage{mathtools}
\usepackage{theoremref}
\usepackage{enumitem}
\usepackage{geometry}
\usepackage{multirow}
\usepackage{hyperref}
\usepackage{faktor}
\usepackage{comment}
\usepackage[new]{old-arrows}
\usetikzlibrary{calc, fit, backgrounds, decorations.markings, arrows.meta}

\makeatletter
\let\@fnsymbol\@alph
\makeatother

\title{Realizing Gruenberg-Kegel graphs of $T$-solvable groups with structurally simplified extensions of $T$}
\author{Lucas Alland\thanks{lalland1@swarthmore.edu; Department of Mathematics, Swarthmore College, 500 College Ave., Swarthmore, PA, 19081, USA.}\and Andrei Fridman\thanks{afridman@hamilton.edu; Department of Mathematics, Hamilton College, 198 College Hill Rd., Clinton, NY, 13323, USA.}\and Thomas Michael Keller\thanks{keller@txstate.edu; Department of Mathematics, Texas State University, 601 University Drive, San Marcos, TX, 78666, USA.}}
\date{}

\DeclareMathOperator{\Aut}{Aut}
\DeclareMathOperator{\pg}{\Gamma}
\DeclareMathOperator{\pgc}{\overline{\Gamma}}
\DeclareMathOperator{\PSL}{PSL}
\DeclareMathOperator{\Sz}{Sz}

\newtheorem{counter}{counter}[subsection]

\theoremstyle{definition}
\newtheorem{definition}[counter]{Definition}

\theoremstyle{plain}
\newtheorem{theorem}[counter]{Theorem}
\newtheorem{lemma}[counter]{Lemma}
\newtheorem{corollary}[counter]{Corollary}
\newtheorem{proposition}[counter]{Proposition}
\newtheorem{conjecture}[counter]{Conjecture}
\newtheorem{theoremmain}{Theorem}
\renewcommand\thetheoremmain{\Alph{theoremmain}}

\theoremstyle{remark}

\begin{document}
\maketitle
\begin{center}
    \vspace{-13pt}
    Dedicated to the memory of I. Martin Isaacs (1940-2025).\vspace{13pt}
\end{center}
\begin{abstract}
   Given a finite group $G$, its prime graph $\Gamma(G)$ (also known as its Gruenberg-Kegel graph) is the graph whose vertices are the prime divisors of $|G|$ and where edges $\{p, q\}$ exist whenever $G$ contains an element of order $pq$. We continue the study of prime graphs for $T$-solvable groups; that is, groups whose composition factors are either abelian or isomorphic to some fixed non-abelian simple group $T$. For a large class of non-abelian simple groups $T$, we prove that the prime graph complements of $T$-solvable groups are always realizable by a solvable group and a quasi simple or almost simple $T$-solvable group acting by automorphisms on a direct product of elementary abelian groups. We conjecture that a similar result holds in full generality. Moreover, we apply our result to classify in purely graph-theoretic terms the prime graph complements of $\PSL(2,13)$-solvable groups, and indicate other interesting classes of groups matching the assumptions of our main theorem.
\end{abstract}

\hfill\break

\hfill\break

\section*{Introduction}
\begingroup
\renewcommand{\thecounter}{\thesection.\arabic{counter}}
If $G$ is a finite group, then let $\pi(G)$ denote the set of prime divisors of $|G|$. The \emph{prime graph} or \emph{Gruenberg-Kegel graph} of $G$, denoted by $\Gamma(G)$, is the graph with vertex set $\pi(G)$, where an edge joins primes $p$ and $q$ if and only if $G$ contains an element of order $pq$. In 2015, \cite{2015REU} achieved a full characterization of the prime graphs of solvable groups, which we cite as a theorem.
\begin{theorem}[\protect{\cite[Theorem 2]{2015REU}}]
    \thlabel{solvableclassification}
    A graph $\Xi$ is isomorphic to the prime graph complement of a solvable group if and only if $\Xi$ is triangle free and $3$-colorable.
\end{theorem}

Following this result, the authors of \cite{2021REU,2022REU,2023REU,Suz,2024PSL} have developed methods for understanding and classifying the prime graphs of larger classes of groups, namely those whose composition factors are either abelian or isomorphic to a fixed non-abelian simple group $T$. We call such groups \emph{$T$-solvable}. The paper \cite{2021REU} classified the prime graphs of $A_5$-solvable groups, \cite{2022REU} classified the prime graphs of all $T$-solvable groups where $|T|$ has exactly $3$ distinct prime divisors, and \cite{2023REU,Suz,2024PSL} classified many $T$-solvable groups where $|T|$ has exactly $4$ distinct prime divisors. These classifications generally take a typical form which we can illustrate by our own classification from Subsection \ref{PSL}.
\begingroup
\renewcommand{\thecounter}{\ref{pslclassification}}
\begin{theorem}
    A graph $\Xi$ is isomorphic to the prime graph complement of some $\PSL(2,13)$-solvable group if and only if one of these conditions holds:
    \begin{enumerate}
        \item $\Xi$ is triangle-free and 3-colorable.
        \item $\Xi$ contains a subset $X$ of four vertices such that in some proper 3-coloring of $\Xi$, the closed neighborhood $N(X) \setminus X$ is monochromatic, and $X$ satisfies one of the following:
        \begin{enumerate}
            \item $X$ contains at most two triangles, all triangles of $\Xi$ are contained within $X$, three vertices of $X$ have no neighbors outside of $X$, and $\Xi[X]$ is connected. Moreover, if a vertex in $X$ indeed has neighbors outside $X$, then it is adjacent to all other vertices in $X$.
            \item There exists a labeling of the vertices in $X$ by $\{a,b,c,d\}$, such that all triangles in $\Xi$ contain the $b-c$ edge, $a$ has no neighbors within $X$, $d$ has no neighbors outside of $X$, $N(c)\setminus\{b\}\subseteq N(b)$, and one of the following is satisfied:
            \begin{enumerate}[label = (\roman*)]
            \item $N(c)\setminus X \subseteq N(a)$
            \item $N(a) = \emptyset$.
            \end{enumerate}
        \end{enumerate}
    \end{enumerate}
\end{theorem}
\endgroup
Since all solvable groups are also $T$-solvable, the interesting cases of these results arise from groups with at least one composition factor isomorphic to $T$. As in the earlier papers, we call such groups \emph{strictly $T$-solvable}. 

In each of \cite{2021REU,2022REU,2023REU,Suz,2024PSL}, the proof of such a classification result proceeds in two directions. First, one must examine the structure of $T$-solvable groups to find a list of necessary conditions satisfied by every corresponding prime graph. Then, one must start with an arbitrary graph satisfying those conditions and construct a $T$-solvable group with an isomorphic prime graph. 

We observe that the proofs of this second step have  always proceeded similarly: begin with either a $T$-solvable quasi simple or almost simple group $H$ whose prime graph is explicitly calculable, and then realize the desired graph by having $H$ and some constructed solvable group both act by automorphisms on a direct product of elementary abelian groups.\footnote{In the case of $A_7$-solvable groups, two non-quasi simple or almost simple groups were used in this step. However, by the methods developed later in \cite{Suz}, this classification could have been done instead in the manner described here.}

We prove that, under certain technical conditions, the above process is always sufficient to realize a given prime graph. Effectively, given these conditions, we reduce the study of prime graphs of $T$-solvable groups to the study of the modular representations of $T$-solvable quasi simple and almost simple groups. In keeping with prior work in these classification results, this paper considers the prime graph complements, denoted by $\pgc(G)$, instead of the prime graphs themselves. Our main result is the following:

\begin{theoremmain}
\thlabel{structuralresult}
    Let $T$ be a non-abelian simple group with Schur multiplier $1$ or $2$ and $\operatorname{Out}(T)$ a—possibly trivial—$2$-group. Let $G$ be a finite strictly $T$-solvable group with at least one edge within $\pi(T)$ in $\pgc(G)$. Suppose that $G$ satisfies at least one of the following:
    \begin{enumerate}
        \item $2$ is adjacent to some $p\in\pi(G)\setminus\pi(T)$ in $\pgc(G)$.
        \item If $t\in\pi(T)$ is Fermat, then for all odd $q \in \pi(G)$ such that $t$ is adjacent to $q$ in $\pgc(G)$, there exists a non-Fermat $r\in\pi(T)$ such that $r$ is adjacent to $q$ in $\pgc(G)$.
    \end{enumerate} Then there exist solvable groups $A$ and $K$ such that $A$ is a direct product of elementary abelian groups, $\pi(K)$ is disjoint from both $\pi(T)$ and its neighbors in $\pgc(G)$, and one of the following is satisfied:
    \begin{enumerate}
        \item A perfect central extension of $T$ by $C_2$ does not appear as a section of $G$, $G$ has a $T$-solvable almost simple quotient group $H$, and $H$ and $K$ act by automorphisms on $A$ such that $\pgc(G)=\pgc(A\rtimes(H\times K))$.
        \item $G$ has a section $E$ which is a perfect central extension of $T$ by $C_2$ and $E$ and $K$ act by automorphisms on $A$ such that $\pgc(G) = \pgc(A\rtimes (E\times K))$.
    \end{enumerate}
\end{theoremmain}

Previous work demonstrates examples which satisfy the results of our main theorem but not all of our hypotheses, such as $\Sz(8)$-solvable groups (the prime graphs of which are classified in \cite{Suz}). We believe these hypotheses can be loosened, and indeed, we conjecture that a similar reduction is always possible. We call a group $H$ \emph{almost quasi simple} if $H$ admits a normal series $1\unlhd Z \unlhd E \unlhd H$ such that $E$ is a quasi simple group and $H/Z$ is an almost simple group. We now conjecture the following.
\begin{conjecture}
\thlabel{generalizedconjecture}
    Let $G$ be a $T$-solvable group. Then either $\pgc(G)$ is triangle free and 3-colorable or there exist solvable groups $A$ and $K$ and a $T$-solvable almost quasi simple group $H$ such that $A$ is a direct product of elementary abelian groups, $\pi(K)$ is disjoint from both $\pi(T)$ and its neighbors in $\pgc(G)$, and $H$ and $K$ act on $A$ such that $\pgc(G)=\pgc(A\rtimes(H\times K))$.
\end{conjecture}

Albeit hard to prove, philosophically speaking, a result like Conjecture \ref{generalizedconjecture} is to be expected to hold: since the map assigning the prime graph to a group is highly non-injective, all prime graphs should be realizable by groups with a limited structure. This became evident in \cite{2015REU}, where it was shown that prime graphs of solvable groups can always be realized by groups of derived length at most three, all of whose Sylow subgroups are elementary abelian. We believe \thref{generalizedconjecture} is the corresponding version for $T$-solvable groups.

Even with the hypotheses of \thref{structuralresult}, our main result simplifies the study of prime graph complements for many interesting $T$-solvable groups. Section \ref{classifications} is dedicated to some such applications. Specifically, we wish to demonstrate that \thref{structuralresult} is highly applicable to the problem of classifying prime graph complements of $T$-solvable groups. While previous work \cite{2022REU,2023REU,Suz,2024PSL} has tackled such problems without a similar structural result, \thref{structuralresult} allows us to pursue more complicated classifications. The following theorem suggests a number of as of yet unclassified finite non-abelian simple groups $T$ for which every $T$-solvable group either has a triangle free and $3$-colorable prime graph complement or satisfies the hypotheses of our main result. The proof may be found in Subsection \ref{applicablegroups}.
\begingroup
\renewcommand{\thecounter}{\ref{thmapplicablegroups}}
\begin{theorem}
    Let $T$ be equal to one of $A_{11}$, $A_{13}$, $A_{14}$, $A_{15}$, $A_{16}$, $\PSL(2,13)$, $S_6(3)$, $HS$, or $G_2(4)$. Then if $G$ is a $T$-solvable group, either $\pgc(G)$ is triangle free and $3$-colorable or $G$ satisfies the hypotheses of \thref{structuralresult}.
\end{theorem}
\endgroup
The prime graphs of $T$-solvable groups for these groups have previously been considered too difficult to classify, either because they have larger and more complicated prime graph complements than considered in previous classifications—that is, they have more than four distinct prime divisors—or, as we will see is true in the case of $\PSL(2,13)$-solvable groups, the associated prime graph complements may contain an arbitrary number of triangles. The only previously identified groups for which this was possible are the $A_5$-solvable ones, whose classification—first carried out in \cite{2021REU} and later redone in \cite[Appendix A]{2023REU}—is significantly more complicated than that of other groups with three prime divisors. However, with the simplifications yielded by \thref{structuralresult}, in Subsection \ref{PSL} we are able to classify the prime graph complements of all $\PSL(2,13)$-solvable groups without much difficulty. With \cite{2024PSL}, this result continues the classification of $\PSL(2,q)$ groups where $|\PSL(2,q)|$ is divisible by exactly $4$ distinct primes.

While classification problems are of particular interest in this paper, we believe \thref{structuralresult} to also be of interest to other questions within the study of prime graphs. We say a group $G$ is \emph{recognizable} by its prime graph if $\pg(G) = \pg(H)$ implies $G\cong H$ for all groups $H$. The authors of \cite{maslova} study various criteria for recognizability, and the following corollary follows immediately from our main result.
\begin{corollary}
    Let $G$ be a group satisfying the hypotheses of \thref{structuralresult} but not of the form $A\rtimes (H\times K)$ or $A\rtimes (E\times K)$, where $A$, $H$, $E$, and $K$ are as given in \thref{structuralresult}. Then $G$ is not recognizable by its prime graph.
\end{corollary}
We also propose that the techniques used in this paper will be of interest to anyone wishing to prove \thref{generalizedconjecture}, and to this end all lemmas are written with generality in mind. The overall strategy employed here of ``gathering'' the relevant representations is likely to prove useful, but more work seems needed in the following areas:
\begin{enumerate}
    \item The technical requirement on $G$ involving Fermat primes comes from the exceptional case of \cite[Theorem A]{Flavell} which we frequently apply by way of \cite[Lemma 2.1.6]{2023REU} to obtain an action of a perfect central extension of $T$ onto a desired chief factor of $G$ (see \thref{flavellapplication,s-action}). This exceptional case indeed occurs, so it seems likely that either the exceptional case needs to be handled to also produce the desired action, or further understanding must be developed as to when the technical requirement is satisfied. The case of $A_5$-solvable groups demonstrates the sort of arguments which may be needed when the exceptional case does occur since the Fermat prime $5\in \pi(A_5)$ often connects to many vertices in the prime graph complements of $A_5$-solvable groups. See \cite{2022REU} and \cite[Appendix A]{2023REU}.
    \item We make frequent use of \cite[Lemmas 2.1.1]{2023REU} to find a normal series for $G$ where an almost simple $T$-solvable group is the first quotient. However, when we do so, sections of $G$ which are perfect central extensions of $T$ may not appear in this normal series. By restricting in this paper to Schur multiplier $2$, we ensure that the only possible non-trivial perfect central extension of $T$ is $2.T$. Furthermore, \thref{2connectsoutward} provides a measure of control over whether $2.T$ appears in our normal series for $G$. To allow for larger Schur multipliers of $T$, it seems one will need stronger methods of this sort to give increased understanding as to how the various possible perfect central extensions of $T$ appear in these normal series.
    \item By restricting $\operatorname{Out}(T)$ to be a $2$-group, this paper sidesteps any issues which may arise when the structure of the prime graph complement is influenced by both the outer automorphisms of $T$ and a non-trivial perfect central extension of $T$.  We do not believe this simplification to be possible in general, so a generalization will likely require using these pieces to build together an almost quasi simple group which realizes the structures of both its constituent sections.
\end{enumerate}

\subsection*{Preliminaries}
All groups in this paper will be assumed to be finite. We will make use of the following notational conventions to align with previous work on this topic.
\begin{enumerate}
    \item For a graph $\Xi$, a vertex set $V$, and an edge set $E$, we write $\Xi \setminus V$ to denote the induced subgraph of $\Xi$ with the vertices from $V$ removed and $\Xi\setminus E$ to denote the subgraph of $\Xi$ with the edges from $E$ removed.
    \item For a graph $\Xi$ and vertex set $V$, we write $\Xi[V]$ to denote the induced subgraph of $\Xi$ on the vertices of $V$.
    \item For two graphs $\Xi_1$ and $\Xi_2$, the notation $\Xi_1\subseteq \Xi_2$ denotes $\Xi_1$ is a subgraph of $\Xi_2$.
    \item For a graph $\Xi$, $v, u$ in the vertex set of $\Xi$, and $\{v,u\}$ in the edge set of $\Xi$, we write $v\in \Xi$ and $v-u\in\Xi$.
    \item For a directed graph $\Vec\Xi$, $v, u$ in the vertex set of $\Vec\Xi$, and a directed edge $(v,u)$ in the edge set of $\Vec\Xi$, we write $v\rightarrow u \in \Vec\Xi$.
    \item For a group $G$ and a set of primes $\mathcal P \subseteq \pi(T)$, the set $N(\mathcal P)$ denotes the vertices adjacent to $\mathcal{P}$ in $\pgc(G)$. If $\mathcal P = \{p\}$, we write $N(p)$ instead of $N(\mathcal P)$.
    \item For a group $G$, $\operatorname{Fit}(G)$ denotes the Fitting subgroup of $G$ and $\Phi(G)$ denotes the Frattini subgroup of $G$.
    \item For groups $H$ and $G$, we write $H \varlonghookrightarrow G$ if $H$ is isomorphic to some subgroup of $G$. 
    \item We often use ATLAS notation to write normal series. If a group $G$ has a normal series, $1=G_0 \unlhd G_1 \unlhd \ldots \unlhd G_k = G$, then we may write $G\cong X_1.X_2.\cdots.X_k$ where each $X_i \cong G_i/{G_{i-1}}$. We also always write a perfect central extension of a finite simple group $T$ by a cyclic group of prime order $C_p$ with $p.T$.
\end{enumerate}
Furthermore, we recall the following definitions from \cite{2015REU}.

\begin{definition}
    A group $G$ is said to be $2$-Frobenius if there exists a normal series $1\unlhd F_1 \unlhd F_2 \unlhd G$ such that $G/F_1$ and $F_2$ are Frobenius groups with kernels $F_2/F_1$ and $F_1$, respectively. We call $G/F_1$ and $F_2$ the upper and lower Frobenius groups.
\end{definition}
\begin{definition}
    A Frobenius group is said to be of type $(p,q)$ if its Frobenius complement is a $p$-group and its kernel is a $q$-group. A $2$-Frobenius group is said to be of type $(p,q,r)$ if the upper Frobenius group is of type $(p,q)$ and the lower Frobenius group is of type $(q,r)$.
\end{definition}
\begin{definition}
    For a solvable group $N$, the Frobenius digraph of $N$, $\vec\Gamma(N)$, is the following orientation of $\pgc(N)$. For each edge $p-q \in \pgc(N)$, the corresponding Hall $\{p,q\}$-subgroup $H\leq N$ is either Frobenius or $2$-Frobenius by \cite[Theorem A]{Williams}. The edge $p-q \in \pgc(N)$ is oriented  $p\rightarrow q$ in $\vec\Gamma(N)$ if and only if $H$ is of type $(p,q)$ or $(p,q,p)$.
\end{definition}

The following lemma is fundamental and will be used throughout this paper without reference.
\begin{lemma}
    Let $G$ be a group and $H$ be a section of $G$. If $p,q\in\pi(H)$ then $p-q\in\pgc(G)$ implies $p-q\in\pgc(H)$.
\end{lemma}

Finally, we clarify that fixed points are always meant to be non-trivial. That is, whenever we mention a group element $g\in G$ \emph{acting with fixed-point} on a group $H$ or a module $V$, we mean there exists $x\in H$ or $x\in V$ such that $x^g = x$ and $x\neq 1_H$ or $0_V$. In case of $G$ acting on a module $V$, we may instead refer to the representation of $G$ and say that $g$ acts with a fixed point \emph{through} the representation.
\endgroup
\section{Reduction of Prime Graph Complements}
In this section, we work towards a proof of \thref{structuralresult}. We begin in Subsection \ref{generallemmas} with some fundamental lemmas which we will use throughout. Subsection \ref{N.T} will consider only $T$-solvable groups which have $T$ as a factor group and will prove \thref{structuralresult} in this reduced case in two parts: considering the edges of $\pgc(G)$ external to $\pi(T)$ and considering those internal to $\pi(T)$. Subsection \ref{induction} will use these partial results to complete the proof of our main theorem.
\subsection{General Lemmas}
\label{generallemmas}
We first consider the structure of $\{r,p\}$-groups with an edge $r-p$ in their prime graph complement. This situation has been well understood in previous work, but it is useful to us to state it as an explicit lemma.
\begin{lemma}
    \thlabel{FrobeniusRef}
    Let $H$ be a $\{r,p\}$-group and $r-p\in\pgc(H)$. Then $H$ is either Frobenius or $2$-Frobenius. Furthermore,
    \begin{enumerate}
        \item If $H$ has a normal $p$-subgroup, then $H$ is Frobenius of type $(r,p)$ or $2$-Frobenius of type $(p,r,p)$.
        \item If there exists $N\unlhd H$ such that $H/N$ is an $r$-group, then $H$ is Frobenius of type $(r,p)$ or $2$-Frobenius of type $(r,p,r)$.
        \item If both (1) and (2) hold, then $H$ is Frobenius of type $(r,p)$.
    \end{enumerate}
\end{lemma}
\begin{proof}
    By Burnside's theorem, $H$ is solvable, so by \cite[Theorem A]{Williams}, $H$ is either Frobenius or $2$-Frobenius. 
    
    For (1), assume $H$ has a normal $p$-subgroup $P$. Now, suppose $H=NK$ is Frobenius, where $N$ is the kernel and $K$ is a complement. As $P$ is nilpotent, we know $P \leq \operatorname{Fit}(H) = N$, so $N$ is not an $r$-group. Since $(|N|, |K|)=1$, we know $H$ is of type $(r,p)$. Suppose instead that $H$ is $2$-Frobenius. By definition, there is a normal series $1\unlhd F_1 \unlhd F_2 \unlhd H$ such that $H/F_1$ and $F_2$ are both Frobenius groups with kernels $F_2/F_1$ and $F_1$, respectively. Since $F_1 = \operatorname{Fit}(F_2)$ and $F_2/F_1 = \operatorname{Fit}(H/F_1)$ we know $F_1 = \operatorname{Fit}(H)$. As $P\unlhd H$ is nilpotent, $P \leq F_1$. Since $F_2/F_1$ acts coprimely on $F_1$, we have that $F_1$ is a $p$-group and $F_2/F_1$ is an $r$-group. This necessitates that $H/F_2$ is a $p$-group, so $H$ is $2$-Frobenius of type $(p,r,p)$.

    For (2), assume $N\unlhd H$ such that $H/N$ is an $r$-group. Suppose $H=N_HK$ is Frobenius, where $N_H$ is the kernel and $K$ is a complement. By \cite[Theorem 8.16]{HuppertEnglish}, either $N\leq N_H$ or $N_H\leq N$. In the first case, as $(|N_H|, |K|)=1$, we have $N= N_H$, so we are done. If $N_H\leq N$, then
    \begin{align*}
        \faktor{(H/N_H)}{(N/N_H)} \cong H/N
    \end{align*}
    is an $r$-group, and hence so is $H/N_H$. Now suppose for a contradiction $H$ is $2$-Frobenius, but of type $(p,r,p)$ and let $1\unlhd F_1\unlhd F_2 \unlhd H$ be the normal series given by definition. Then $F_1$ is a $p$-group, which implies $F_1\unlhd N$. So $N/F_1 \unlhd H/F_1$ and by \cite[Theorem 8.16]{HuppertEnglish}, either $N/F_1 \leq F_2/F_1$ or $F_2/F_1 \leq N/F_1$, both of which lead to a contradiction by the argument above.

    Finally, (3) follows immediately by noticing that $H$ cannot be both $(r,p,r)$ and $(p,r,p)$.
\end{proof}

In a similar style to \cite[Subsection 2.2]{2023REU}, we establish the following two similar criteria for disallowing edges $r-p$ where $r\in\pi(T)$ and $p\in\pi(G)\setminus \pi(\Aut(T))$. While we will often use these next two lemmas in tandem, their proofs differ almost entirely.
\begin{lemma}
\thlabel{solvable-rprime}
    Let $T$ be a non-abelian simple group, $G$ be strictly $T$-solvable, and $K\leq G$ be of the form $K\cong N.T$ as given by \cite[Lemma 5.2]{Suz}. For all $r\in\pi(T)$ such that $r$ does not divide the Schur multiplier of $T$, if $r-p \in \pgc(G)$ for some $p\in\pi(G)\setminus\pi(\Aut(T))$, then $r\nmid|N|$
\end{lemma}
\begin{proof}
    Suppose $r \mid |N|$, but $r - p \in \pgc(G)$ for some $p \not\in\pi(T)$. It will be sufficient to find our contradiction in $K$, so without loss of generality assume $G=K$. Let $L \leq G$ be the solvable subgroup given by $L\cong N.R$ where $R$ is a Sylow $r$-subgroup of $T$. We see $r-p\in\pgc(L)$, so by \thref{FrobeniusRef}, $H \in \operatorname{Hall}_{\{r,p\}}(L)$ is Frobenius of type $(r,p)$ or 2-Frobenius of type $(r,p,r)$. In the first case, the Sylow $r$-subgroups of $G$ must be cyclic or generalized quaternion. In the latter case, by \cite[Corollary 2.3]{2015REU}, the Sylow $p$-subgroups of $G$ are cyclic. We consider each of these individually. 

    {\bf Case 1a: The Sylow $r$-subgroups of $G$ are cyclic.} Take a chief series for $G$ given by $1=N_{0}\unlhd N_1\unlhd \ldots \unlhd N_k=N\unlhd G$
    where $\ell$ is the largest integer such that $N_\ell/N_{\ell-1}$ is a $r$-group. Since it is enough to show $r-p\not\in\pgc(G/N_{\ell-1})$, we may assume $N_{\ell-1}=1$. Now, $N_\ell\cong C_{r}$ and $G/N_\ell$ acts by conjugation on $N_\ell$. Let $W/N_\ell \unlhd G/N_\ell$ be the kernel of this action. Then $(G/N_\ell)/(W/N_\ell) \cong G/W \varlonghookrightarrow \Aut(N_\ell) \cong C_{r-1}$. So $W$ contains all the order $r$ elements from $G$. Thus, $W\not\leq N$, so $W/(W\cap N) \neq 1$. Furthermore, since 
    \begin{align*}
        W/(W\cap N) \cong WN/N \unlhd G/N \cong T,
    \end{align*}
    it follows that $W/(W\cap N) \cong T$. Let $M = W\cap N$. Define $Q =  M/N_\ell$. Since $N_\ell$ contains all elements of order $r$ in $M$, $Q$ contains no element of order $r$. By the Schur-Zassenhaus theorem, $M \cong {N_\ell} \rtimes Q$. Since $Q$ is contained within the kernel of the action of $K/N_\ell$ on $N_\ell$, we know the action from $Q$ is trivial. That is, $M\cong N_\ell\times Q$. Next, observe that $(W/Q)/(M/Q) \cong W/M \cong T$, and $M/Q \cong C_{r}$. Moreover,
    \begin{align*}
        \faktor{\mathbf{N}_{W/ Q}( M/ Q)}{\mathbf C_{ W/ Q}( M/ Q)} = \faktor{( W/ Q)}{\mathbf C_{ W/ Q}( M/ Q)} \varlonghookrightarrow  \Aut( M/ Q) \cong C_{r-1}.
    \end{align*}
    Since $ M/ Q \leq \mathbf C_{ W/ Q}( M/ Q)$, $\mathbf C_{ W/ Q}( M/ Q)$ is normal in $W/ Q$, and $ M/ Q$ is a maximal normal subgroup of $W/Q$, we have that either $\mathbf C_{ W/ Q}( M/ Q) =  M/ Q$ or $\mathbf C_{ W/ Q}( M/ Q) =  W/ Q$. The first case is impossible, as $( W/ Q)/( M/ Q) \cong T$ is not isomorphic to a subgroup of  $C_{r-1}$. Hence $ M/ Q = Z( W/ Q)$, so $ W/ Q$ is a central extension of $T$. 

    Suppose $ W/ Q$ is a split extension. Then since $T$ is simple, its action on $M/Q$ is trivial. So we have $ W/ Q \cong  M/ Q \times T$, which is a contradiction as the Sylow $r$-subgroups of $W/Q$ are cyclic. It follows that $ W/ Q$ is a non-split central extension of $T$ by $ M/ Q \cong C_{r}$. Then the commutator subgroup, $(W/Q)^\prime$, is a perfect central extension of $T$ by $C_{r}$, which contradicts that $r$ does not divide the Schur multiplier of $T$. This completes the first case.

    {\bf Case 1b: The Sylow $r$-subgroups of $G$ are generalized quaternion.} First, we remark that $r=2$. Let $Q$ be a Sylow $2$-subgroup of $G$. We see that $Q/(Q\cap N)$ is a Sylow $2$-subgroup of $T$. By Cayley's normal $2$-complement theorem, $Q/(Q\cap N)$ is not cyclic. Thus, since furthermore $Q\cap N$ is non-trivial, $|Q:Q\cap N|\geq 4$ which implies $Q\cap N$ is cyclic. Now, if we take a chief series for $G$ through $N$, every characteristic $2$ chief factor must be isomorphic to $C_2$—which reduces back to Case 1a. 

    {\bf Case 2: The Sylow $p$-subgroups of $G$ are cyclic.} Take the chief series for $G$ given by
    $1=N_0 \unlhd N_1\unlhd \ldots \unlhd N_k=N\unlhd G$,
    and let $\ell$ be the maximal integer such that $N_\ell/N_{\ell-1}$ is a $p$-group. Again, it suffices to show that $r-p\not\in\pgc(G/N_{\ell-1})$, so we assume $N_{\ell-1}=1$. Now, ${N_\ell} \cong C_p$ is abelian, so $ G/ N_\ell$ acts by conjugation on ${N_\ell}$. Let $ W/ N_\ell$ be the kernel of this action. If $ W\leq  N$, then by
    \begin{align*}
        T\cong  G/ N\cong\faktor{(G/ W)}{( N/ W)},
    \end{align*}
    $T$ is a factor of $G/W$. Recall $G/W\varlonghookrightarrow \operatorname{Aut}( N_\ell)$ is abelian, so $T$ is abelian—a contradiction. Thus, suppose $W$ is not a subgroup of $N$. Since 
    \begin{align*}
        1\neq W/( W\cap N) \cong  W N/ N \unlhd  G/ N \cong T,
    \end{align*}
    we have that 
    \begin{align*}
        \faktor{(W/{N_\ell)}}{(( W\cap {N})/N_\ell)}\cong  W/( W\cap  N) \cong T.
    \end{align*}
    Then the kernel $ W/ N_\ell$ contains an order $r$ element which acts trivially on the $p$-group $N_\ell$, a contradiction with $r-p\in\pgc(G)$ by \cite[Corollary 1.5]{Suz}.
\end{proof}

The next lemma provides a similar result when $r=2$ does divide the Schur multiplier, but also provides structural control over the group as a whole.
\begin{lemma}
\thlabel{2connectsoutward}
    Let $T$ be a non-abelian simple group with a perfect central extension $2.T$. Suppose $G$ is a strictly $T$-solvable group of the form $G\cong N.T.O$ for $N$ solvable and $O\leq\operatorname{Out}(T)$. Let $2-p \in \pgc(G)$ for some $p\in\pi(G)\setminus\pi(T)$. Then $G\cong L.2.T.O$ where $L$ is a solvable $2^\prime$ group and $2.T$ is a perfect central extension and $2-r\not\in\pgc(G)$ for all $r\in\pi(T)$.
\end{lemma}
\begin{proof}
    The proof is a variation of 
    \cite[Lemma 1.1]{2024PSL}. At the start of the second paragraph where the authors assume $G \cong N.T$ by passing to the subgroup $K$ given by \cite[Lemma 2.1.2]{2023REU}, we instead work in $N.T\unlhd G$ as assumed exists by the assumption $G\cong N.T.O$. The proof then proceeds identically to show that $N.T\cong L.2.T$. Then, still following the proof of \cite[Lemma 1.1]{2024PSL}, since we have $L\operatorname{char} N$ and $C_2\operatorname{char} N.T/L$, we know $L\unlhd G$ and $C_2\unlhd G/L$. Thus, $G\cong L.2.T.O$.
\end{proof}

What follows is a similar edge-nonexistence criteria to the proceeding two lemmas, except we consider edges between elements of $\pi(T)$.
\begin{lemma}
\thlabel{divisibilityadjacency}
    Let $T$ be a non-abelian simple group, $G$ be strictly $T$-solvable, and $K\leq G$ be of the form $K\cong N.T$ as given by \cite[Lemma 5.2]{Suz}. Suppose that $s, r \in \pi(T)$ both divide the order of $N$. Then $r-s\not\in\pgc(G)$.
\end{lemma}
\begin{proof}
        Suppose for the sake of contradiction $r-s\in\pgc(G)$. It suffices to find a contradiction in $K\leq G$, so without loss of generality suppose $K=G$. Take a chief series for $G$: $1=N_0\unlhd N_1 \unlhd \ldots \unlhd N_k=N \unlhd G$. Let $\ell$ be the maximal integer such that $r,s \mid|N_{\ell-1}|$. Then, $N_\ell/N_{\ell-1}$ is either an $r$ or $s$-group. Without loss of generality, assume the latter. Since $G/N_{\ell-1}$ satisfies the hypotheses of the claim and is sufficient to form our contradiction, without loss of generality let $N_{\ell-1} = 1$. Suppose $r$ is odd. Then we reach a contradiction by applying \cite[Lemma 2.1.3]{2023REU} with their $p_i$ replaced by our $r$, and their $s$ and $p_j$ replaced by our $s$. So we know $r=2$.
        
        Let $R$ be the Sylow $2$-subgroup of $T$. Then $N.R$ is a solvable group, so its Hall $\{2,s\}$-subgroup $H$ is Frobenius of type $(2,s)$, by \thref{FrobeniusRef}. We have $H/(H\cap N) \cong R$ and, by Cayley's normal $2$-complement theorem, that $R$ is not cyclic. Thus, write $H=SQ$ where $S$ is an $s$-group and $Q$ is a generalized quaternion group. Let $Q_0 = Q\cap N$. We see $Q_0$ is non-trivial since $2\mid|N|$. Then, since $Q/Q_0$ is not cyclic, $|Q:Q_0|\geq 4$, which implies $Q_0$ is cyclic. By Cayley's normal $2$-complement theorem, we may write $N \cong L.Q_0$ for $L\unlhd N$ of odd order. Since $L$ is characteristic in $N$, it is normal in $G$. Thus, $E\cong(G/L)/(N/L)$ acts on $N/L \cong Q_0$ by conjugation. The kernel of this action must contain all the order $s$ elements of $E$ since $\Aut(Q_0)$ is a $2$-group—which contradicts that $2-s\in\pgc(G)$. 
\end{proof}

We now wish for a generalization of \cite[Lemma 2.3.6]{2023REU} which provides equality of labeled graphs rather than unlabeled isomorphism. We first need the following result:
\begin{lemma}
\thlabel{frobeniusexistence}
    Let $p$ and $q$ be distinct primes. Then there exists an elementary abelian $q$-group $Q$ and an action by automorphisms of $C_p$ on $Q$ such that $Q\rtimes C_p$ is Frobenius.
\end{lemma}
\begin{proof}
    By Fermat's little theorem, $p$ divides $q^{p-1}-1$. Let $F$ be a finite field of order $q^{p-1}$ and consider the additive group $F^+$. We see the multiplicative group $F^\times$ is of order $q^{p-1}-1$ and acts Frobeniusly on $F^+$. Then there is an order $p$-subgroup of $F^\times$ acting Frobeniusly on $F^+$. We identify this order $p$-subgroup with $C_p$ and let $Q=F^+$. By construction $Q\rtimes C_p$ is Frobenius, and $Q$ is an elementary abelian $q$-group.
\end{proof}

Using the lemma above, the proof of \cite[Lemma 2.3.6]{2023REU} generalizes without a significant change in the proof strategy. While the general structure is the same, we find it easier to state our result in terms of modular representations rather than complex ones.
\begin{lemma}
\thlabel{solvablestructure}
    Let $\Xi$ be a finite graph with vertices from the set of all primes and $E$ be a group with $\pgc(E)\subseteq \Xi$. Suppose the following hold:
    \begin{enumerate}
        \item $\Xi\setminus\pi(E)$ is triangle free.
        \item $\Xi\setminus\pi(E)$ has a 3-coloring $\{\mathcal O, \mathcal D, \mathcal I\}$ such that vertices in $N(\pi(E))\setminus \pi(E)$ are colored $\mathcal I$.
        \item Given $p\in \Xi\setminus\pi(E)$ there exists a finite module of characteristic $p$ such that every element of order $r\in\pi(E)$ acts Frobeniusly if and only if $r-p\in\Xi$.
    \end{enumerate}
    Then there exists solvable groups $J$ and $K$ such that $J$ is a direct product of elementary abelian groups, and $E$ and $K$ act by automorphisms on $J$ with $\Xi = \pgc(J\rtimes(E\times K))$.
\end{lemma}
\begin{proof}
    We define an orientation on $\Xi\setminus\pgc(E)$ by its three-coloring: $\mathcal O\to\mathcal D$, $\mathcal O \to \mathcal I$, and $\mathcal D\to\mathcal I$. Notice this orientation admits no directed $3$-paths. Since we also recall $\Xi\setminus\pgc(E)$ is triangle free, the orientation thus admits no cycle.

    Partition $\Xi\setminus \pgc(E)$ into sets $\{\mathcal V,\mathcal U,\mathcal W\}$ where vertices in $\mathcal V$ have out-degree $0$, vertices in $\mathcal U$ have directed edges into $\mathcal V$ and vertices in $\mathcal W$ have directed edges into $\mathcal U$. Since our orientation on $\Xi\setminus \pgc(E)$ is acyclic, these subsets of the vertices indeed constitute a partition. We also see $N(\pi(E))\setminus \pi(E) \subseteq \mathcal I \subseteq V$.

    Fix some $u\in \mathcal U$. For each $w\in W$ such that the edge $w\to u$ is in our orientation, we have by \thref{frobeniusexistence} an action of $C_w$ on some $(C_u)^{m_w}$. Let $n_u$ be the least common multiple all the integers $m_w$ for $w\in \mathcal W$ adjacent to $u$. Then we see that each such $C_w$ acts Frobeniusly on the whole group $(C_u)^{n_u}$ by acting Frobeniusly on each of its components. Let $W=\bigtimes_{w\in \mathcal W} C_w$. We define an action by automorphisms of $W$ onto $(C_u)^{n_u}$ by letting each component $C_w$ act Frobeniusly if $w\to u$ is in our orientation—and trivially otherwise. Then, if $U=\bigtimes_{u\in \mathcal U} (C_u)^{n_u}$, we may have $W$ act on $U$ as it does on each of its components.

    By the same argument, for each $v\in \mathcal V$ we have some $n_v$ such that $C_u$ has a Frobenius action on $(C_v)^{n_v}$ for every $u\in\mathcal U$ such that $u\to v$ is in our orientation. Define an action of $U$ on $(C_v)^{n_v}$ by letting each component $C_u$ act Frobeniusly if $u\to v$ is in our orientation—and trivially otherwise. Then if $V=\bigtimes_{v\in\mathcal V}(C_v)^{n_v}$, we may have $U$ act on $V$ as it does on each of its components. 
    
    Let $K = U\rtimes W$ act on $V$ by the action of $U$ on $V$ and the trivial action of $W$. We see by construction that $\Xi\setminus \pgc(E) = \pgc(V\rtimes K)$. Now, take $v\in\mathcal V$ such that $v\in N(\pi(E))$. Let $(C_v)^{l_1}$ be the finite characteristic $v$ module of $E$ guaranteed by hypothesis and let $(C_v)^{l_2}$ be the Sylow-$v$ subgroup of $V$. Then, if $k_v$ is the least common multiple of $l_1$ and $l_2$ we have an action by automorphisms of $E$ and $K$ on $(C_v)^{k_v}$ by their actions on its components. Then, we have an action by automorphisms $E$ and $K$ on the product of elementary abelian groups $J=\bigtimes_{v\in\mathcal V} (C_v)^{k_v}$ such that $\Xi = \pgc(J\rtimes(E\times K))$.
\end{proof}

To align with the statement of \cite[Lemma 2.3.6]{2023REU} we also present the following corollary to the above lemma. This corollary may be of use to anyone considering \emph{labeled} classifications of the prime graph complements of $T$-solvable groups.
\begin{corollary}
    Let $E$ and $\Xi$ match the hypotheses of \thref{solvablestructure} except suppose all the assumed representations of finite characteristic are instead complex. Then the conclusions of the lemma are also satisfied for this $E$ and $\Xi$.
\end{corollary}
\begin{proof}
    The proof is immediate from \cite[Lemma 3.4]{2021REU}, since each of the desired complex representations gives the representation of finite characteristic assumed in \thref{solvablestructure}.
\end{proof}

We close our section on general lemmas with two small facts about fixed points of representations which will be useful to us later.
\begin{lemma}
\thlabel{fixedpointunions}
   Let $E$ be a group and $\mathcal P,\mathcal Q\subseteq\pi(E)$. Suppose $E$ has representations $\rho_{\mathcal P}$ and $\rho_{\mathcal Q}$ over an arbitrary field $F$ where elements of orders from $\mathcal P$ or $\mathcal Q$ respectively act with fixed points and elements of other prime orders act Frobeniusly. Then $\rho_{\mathcal P} \oplus \rho_{\mathcal Q}$ is a $F$-representation of $E$ where elements of orders $\mathcal P\cup \mathcal Q$ act with fixed points and all others act Frobeniusly.
\end{lemma}
\begin{proof}
    Immediate.
\end{proof}

\begin{lemma}
    \thlabel{fixedpointextension}
    Let $E$ be a group and $\mathcal P \subseteq \pi(E)$. Suppose $E$ has a representation $\rho$ over an arbitrary field $F$ where elements of orders from $\mathcal P$ act with fixed points and elements of other prime orders act Frobeniusly. Let $\hat E$ be a group with a normal subgroup $N\unlhd \hat E$ such that $\hat E/N \cong E$. Then $\hat E$ has an $F$-representation onto the image of $\rho$ in which elements of orders from $\mathcal P \cup \pi(N)$ act with fixed points and elements of other prime orders act Frobeniusly.
\end{lemma}
\begin{proof}
    Let $\hat \rho$ be the composition of the projection map from $\hat E$ onto $E$ with the representation $\rho$. All elements of prime orders which do not divide $|N|$ act Frobeniusly through $\hat \rho$ if and only if they do so in $\rho$, and elements of all order which divide $|N|$ appear in the kernel of $\hat \rho$—and hence act with fixed points. 
\end{proof}
\subsection{Groups of the form $N.T$}
\label{N.T}
This subsection effectively proves \thref{structuralresult} for strictly $T$-solvable groups $G$ which have a solvable normal subgroup $N$ such that $G/N \cong T$. We will prove \thref{structuralresult} in two parts, \thref{externalstructure} is essentially the proof for the edges within the solvable section of the prime graph complement and between the solvable and non solvable sections. \thref{internalstructure} separately handles the edges within the non solvable sections. We save the combination of these results for the full proof of \thref{structuralresult}. The lemmas in this section also all deal with $T$-solvable groups of this specific form, although we sometimes wish to allow our upper quotient to be isomorphic to a general perfect central extension of $T$, instead of strictly $T$ itself.

We begin with a generalization of \cite[Lemma 2.3.5]{2023REU} to handle when many primes in $\pi(T)$ are adjacent to primes in $\pi(G)\setminus\pi(T)$. Effectively, we show that the class of groups considered in this subsection all have the prime graph complements with the coloring assumed by \thref{solvablestructure}.
\begin{lemma}
\thlabel{monochromaticneighbors}
    Let $T$ be a non abelian simple group and let $G$ be a group of the form $G\cong N.T$ for $N$ solvable. Let $\Gamma$ be the subgraph of $\pgc(G)$ with all the edges internal to $\pi(T)$ removed. Then $\Gamma$ is $3$-colorable such that $N(\pi(T))$ is monochromatic.
\end{lemma}
\begin{proof}
    This proof is inspired by the proof of the Gallai-Roy-Vitaver Theorem (\cite[Theorem 8.5]{GraphTheory}).

    Let $n = |\pi(T)|$ and $V = \{v_1, v_2, \ldots, v_n\}$ be a set of $n$ abstract vertices. We will construct the digraph $\vec\Gamma$ on the vertex set $\pi(G)\cup V$ as follows. For each edge of $\Gamma$ within $\pi(G)\setminus\pi(T)$, add it to $\vec \Gamma$ oriented by the Frobenius digraph on $N$. For each edge of $\Gamma$ between $\pi(T)$ and $\pi(G)\setminus\pi(T)$, add it to $\vec\Gamma$ oriented from the former to the latter. Enumerate the elements of $\pi(T)$ by $\pi(T) = \{p_1, p_2, \ldots, p_n\}$, and add the oriented edge $v_i\rightarrow p_i$ to $\vec\Gamma$ for each $1\leq i\leq n$.
    
    By \cite[Lemma 2.3.4]{2023REU} and \cite[Corollary 2.7]{2015REU}, $\vec\Gamma$ contains no directed paths of length $3$. Furthermore, since we see $\vec \Gamma$ has the form of isolated paths connected onto a subset of the Frobenius digraph of $N$, and by \thref{solvableclassification} we have that the Frobenius digraph of $N$ is triangle-free, we must have that $\vec \Gamma$ is triangle free. Thus, $\vec\Gamma$ is acyclic. Color $\vec\Gamma$ by the length of the longest directed path terminating in each vertex. Since the graph is acyclic this clearly yields a proper $3$-coloring of $\vec \Gamma$. Furthermore, we have constructed $\vec\Gamma$ so each neighbor of $\pi(T)$ in $\Gamma$ is at the end of a path of length 2. Thus, this construction yields the desired $3$-coloring on our original graph $\Gamma$.
\end{proof}
\begin{figure}[h!]
    \centering
    \begin{tikzpicture}[
      node/.style={circle, draw, minimum size=8mm},
      mid arrow/.style={
        decoration={markings, mark=at position 0.5 with {\arrow[scale=1.2]{>}}},
        postaction={decorate},
        line width=1.2pt,
        >=Stealth
      }]

        \node[node] (B) at (0.5,-0.5) {$p_3$};
        \node[node] (D) at (0.5,2.5) {$p_1$};
        \node[node] (A) at (2,1) {$p_2$};
        
        \node[node] (E) at (3,5) {$v_1$};
        \node[node] (F) at (3,-3) {$v_3$};
        \node[node] (C) at (-1,1) {$v_2$};
        
        \node[node,fill =blue!30] (G) at (5.5, 3.5) {};
        \node[node,fill =blue!30] (H) at (5.5,-1.5) {};
        \node[node,fill =blue!30] (I) at (5.5,1) {};
        \node[node,fill =green!30] (J) at (7.5,2.5) {};
        \node[node, fill =red!30] (K) at (7.5,-1.5) {};
        \node[node, fill =red!30] (L) at (7.5,1) {};
        \node[node, fill =green!30] (M) at (9.5,-0.2) {};
        
        \draw[mid arrow] (C) -- (A);
        \draw[mid arrow] (E) -- (D);
        \draw[mid arrow] (F) -- (B);
        \draw[mid arrow] (D) -- (G);
        \draw[mid arrow] (D) -- (I);
        \draw[mid arrow] (B) to[out=0,in=260] (G);
        \draw[mid arrow] (B) -- (H);
        \draw[mid arrow] (J) -- (G);
        \draw[mid arrow] (L) -- (H);
        \draw[mid arrow] (L) -- (I);
        \draw[mid arrow] (K) -- (H);
        \draw[mid arrow] (M) -- (K);
        \draw[mid arrow] (M) -- (L);
        \draw[mid arrow] (A) -- (I);
        
        \begin{pgfonlayer}{background}
          \coordinate (extendleft) at ($(G) + (-0.6, 0)$);;
          
          \node[
            draw=purple,
            dashed,
            thick,
            rounded corners,
            inner sep=6pt,
            fit=(extendleft)(G)(H)(I)(J)(K)(L)(M),
            name=regionbox
          ] {};
        
          \node[
            anchor=north east,
            rectangle,
            draw=none,
            fill=none,
            font=\small\bfseries,
            text=purple,
            xshift=-4pt,
            yshift=-4pt
          ] at (regionbox.north east) {$\pi(G)\setminus \pi(T)$};
        
            {};
        \end{pgfonlayer}
    \end{tikzpicture}
    \caption{An example of $\vec\Gamma$ as given in \thref{monochromaticneighbors} }
\end{figure}

The following lemma allows us to take effectively construct a characteristic $p>0$ representation of a perfect central extension of $T$ from many such representations of groups $N.E$ which have a perfect central extension $E$ of $T$ as a factor group by a normal subgroup $N$ of $p^\prime$ order. This lemma is an application of the major result \cite[Theorem A]{Flavell} by way of \cite[Lemma 2.1.6]{2023REU} and the assumptions present there and thus in the following necessitate many of the technical requirements on \thref{structuralresult}.
\begin{lemma}
    \thlabel{flavellapplication}
    Let $T$ be a non-abelian simple group with $r\in\pi(T)$ and $G$ be a group with $G\cong P \rtimes (N.E)$ such that
    \begin{enumerate}
        \item Either $r$ is not Fermat or $N/(N\cap \mathbf C_{N.E}(P))$ is a $2^\prime$-group.
        \item $P$ is an elementary abelian $p$-group for some prime $p\neq r$ where $r-p\in\pgc(G)$,
        \item $N$ is an $\{r,p\}^\prime$ group, and
        \item $E$ is a perfect central extension of $T$.
    \end{enumerate}
    Then there exists some $\hat E$ which acts by automorphisms on $P$ satisfying all of the following.
    \begin{enumerate}
        \item $\hat E$ is section of $N.E$, $\hat E$ is a perfect central extension of $T$, $\hat E$ and contains $E$ as a factor group,
        \item $s-p\in\pgc(G)$ implies $s-p\in\pgc(P\rtimes \hat E)$ for all $s\in\pi(T)$.
        \item $s-p\in\pgc(P\rtimes \hat E)$ implies $s-p\in\pgc(G)$ for all $s\in\pi(T)$ such that $s\nmid|N|$.
     \end{enumerate}
\end{lemma}
\begin{proof}
    Let $H = N.E$, let $K = \mathbf C_{H}(P)$, and let $R$ be any subgroup of $H$ isomorphic to $C_r$. We know $r-p \in \pgc(G)$ so $r\nmid|K|$. Thus, by Schur-Zassenhaus,
    \begin{align*}
        RN/(RN\cap K) = RN/(N\cap K) \cong N/(N\cap K) \rtimes R.
    \end{align*}
    We recall by assumption that $R \cong C_r$ and that $R$ does not have a fixed point on $P$. Furthermore, we see that $RN/(RN\cap K)$ acts faithfully on $P$. So, by the contrapositive to \cite[Lemma 2.1.6]{2023REU}, $[R, N/(N\cap K)] = 1$. For all $s \in \pi(T)$ where $s\nmid|N|$, $s-p\in\pgc(G)$ only if elements of orders $s$ in $H/(N\cap K)$ act without fixed points on $P$. Similarly, for arbitrary $s\in\pi(T)$, if elements of orders $s$ in $H/(N\cap K)$ act with fixed points on $P$ then $s-p\not\in\pgc(G)$. So we may assume without loss of generality that $N\cap K = 1$.

    Now define $C = \mathbf{C}_H(N)$. By the above argument, $C$ contains all order $r$ elements of $H$. Since $r\nmid|N|$ and
    \begin{align*}
        C/(C\cap N) \cong CN/N \unlhd H/N \cong E,
    \end{align*}
    we know $C/(C\cap N)$ must be isomorphic to a normal subgroup of $E$ containing all of its order $r$ elements. The only such subgroup is $E$ itself. By the correspondence theorem, let $A\unlhd C$ be such that $A/(C\cap N) \cong Z(E)$. By identifying the Sylow $r$-subgroups of $H$ with those of $E$, they act trivially on both $C\cap N$—since $C\cap N \leq Z(C)$—and on $A/(C\cap N)$—since $A/(C\cap N) = Z(E)$. If $Z(E)$ is an $r^\prime$-group then by \cite[Exercise 3E.3]{IsaacsFiniteGroups}, the Sylow $r$-subgroups of $H$ act trivially on $A$. If $Z(E)$ is an $r$-group then by Schur-Zassenhaus, $A\cong C\cap N \times Z(E)$ and the Sylow $r$-subgroups also act trivially on $A$. In either case, $\mathbf C_C(A) = C$, so $A\leq Z(C)$. If $Z(E)$ is neither an $r^\prime$-group nor an $r$-group then by applying the above argument to its $r$ and $r^\prime$ subgroups we also conclude $A\leq Z(C)$. Since $T$ is simple, we must have $A = Z(C)$. So $C$ is a central extension of $T$. Its commutator subgroup, $C^\prime$ is thus a perfect central extension of $T$. Since $E$ is perfect, we know
    \begin{align*}
        E= E^\prime \cong (C/(C\cap N))^\prime = C^\prime ( C \cap N)/( C \cap N) \cong  C^\prime/( C^\prime\cap N),
    \end{align*}
    and thus $C^\prime$ is the desired group $\hat E$ with the action on $P$ inherited from $H$.
\end{proof}

For each given non-edge between $\pi(T)$ and $\pi(G)\setminus\pi(T)$ in the prime graph complement for a group, we wish to be able to find an action by automorphisms matching the assumptions of \thref{flavellapplication} which has the corresponding fixed point. The following lemma allows us to do so at the cost of passing to a subgroup.
\begin{lemma}
    \thlabel{p-maximal}
    Let $T$ be a non-abelian simple group, $G$ be of the form $G\cong N.E$ with $N$ solvable and $E$ a perfect central extension of $T$, and $p\in\pi(G)\setminus\pi(T)$. For every $s\in\pi(T)$ such that $s-p\not\in\pgc(G)$ and $s\nmid|N|$, there exists a subgroup $H\leq G$ with a chief series $1=H_0 \unlhd H_1 \unlhd \ldots \unlhd H_k \unlhd \ldots \unlhd H$ satisfying
    \begin{enumerate}
        \item $H/H_k \cong E$ and
        \item If $\ell$ is the maximal integer such that $H_\ell/H_{\ell-1}$ is a $p$-group, then $s-p\not\in\pgc(H/N_{\ell-1})$.
    \end{enumerate}
\end{lemma}

\begin{proof}
    Suppose for the sake of contradiction that $G$ is a counterexample of minimal order for some fixed $s\in\pi(T)$ such that $s\nmid|N|$. Take a chief series for $G$ through $N$: $1=N_0\unlhd N_1 \unlhd \ldots\unlhd N_k = N \unlhd \ldots \unlhd G$ and let $\ell$ be the maximal integer where $N_\ell/N_{\ell-1}$ is a $p$-group. Since $G$ is a minimal counterexample, $s-p\in\pgc(G/N_{\ell-1})$. Applying Schur-Zassenhaus, we know $G/N_{\ell-1} \cong N_{\ell}/N_{\ell-1} \rtimes K/N_{\ell-1}$ for some $K\leq G$. Since $s-p\not\in\pgc(G)$, there exist commuting $x,y\in G$ of orders $s$ and $p$, respectively. Since $s-p\in\pgc(G/N_{\ell-1})$, we see that $y \in N_{\ell-1}$. Furthermore, there exists some $g\in G$ such that $g^{-1}xg \in K$. As $N_{\ell-1}\unlhd K$, we know $g^{-1}xg \in K$ commutes with $g^{-1}yg \in K$. So $s-p\not\in \pgc(K)$. We see $K\cong (N\cap K).E$, so by our minimality assumption $K$ has the desired subgroup $H$. Since $H$ is also a subgroup of $G$, we have the desired contradiction.
\end{proof}

The following allows us to find representations of perfect central extensions of $T$ with positive characteristic dividing both $\pi(T)$ and $\pi(G)\setminus\pi(T)$. The associated semidirect product may have more edges than $\pgc(G)[\pi(T)]$ but not fewer. Effectively, when \thref{divisibilityadjacency} guarantees the non-existence of an edge in $G$, we construct a module which will allow us to also guarantee by \thref{divisibilityadjacency} that this same edge is not present in the group we eventually construct.
\begin{lemma}
    \thlabel{removingneighbors}
    Let $T$ be a non-abelian simple group, $G$ be of the form $G\cong N.E$ with $N$ solvable and $E$ a perfect central extension of $T$. Let $s\in\pi(T)$ divide $|N|$. Suppose that furthermore $G$ satisfies at least one of the following:
    \begin{enumerate}
        \item $s=2$.
        \item $N$ is a $2^\prime$ group.
        \item If $t\in\pi(T)$ is Fermat and $t-s\in\pgc(G)$, there exists a non-Fermat $r\in\pi(T)$ such that $r-s\in\pgc(G)$.
    \end{enumerate}
    Then there exists an elementary abelian $s$-group $A$, a section $\hat E$ of $G$ which is a perfect central extension of $T$ containing $E$ as a section, and an action of $\hat E$ on $A$ such that $\pgc(G)[\pi(T)] \subseteq \pgc(A\rtimes \hat E)$.
\end{lemma}
\begin{proof}
    If $N(s)\cap\pi(T)=\emptyset$ then $A=C_s$ with $E$ acting on $A$ trivially suffices. So let $r\in\pi(T)$ be such that $r-s\in\pgc(G)$. Without loss of generality, we may assume that either $s=2$, $N$ is a $2^\prime$-group, or $r$ is not Fermat. Furthermore, by \thref{divisibilityadjacency}, $r\nmid|N|$. Take a chief series for $G$: $1=N_0\unlhd N_1 \unlhd \ldots \unlhd N_k=N \unlhd \ldots \unlhd G$ and let $\ell$ be the maximal integer such that $N_\ell/N_{\ell-1}$ is an $s$-group. Since $\pgc(G)[\pi(T)]\subseteq\pgc(G/N_{\ell-1})[\pi(T)]$, we may assume without loss of generality that $N_{\ell-1} = 1$. Furthermore,  $\pgc(G)[\pi(T)]\subseteq\pgc(N_\ell \rtimes G/N_{\ell})[\pi(T)]$ by \cite[Corollary 1.6]{Suz}, so  assume without loss of generality that $G\cong N_\ell \rtimes G/N_\ell$. Since $N/N_{\ell}$ is an $s^\prime$ group, we may now assume that either $r$ is non-Fermat or $N/N_{\ell}$ is $2^\prime$. Then apply \thref{flavellapplication} to obtain an action of a perfect central extension $\hat E$ of $T$ containing $E$ as a section on $N_\ell$ where for all $t\in\pi(T)$ such that $t-s\in\pgc(G)$ we have $t-s\in\pgc(N_\ell\rtimes \hat E)$.
    Thus, $\pgc(G)[\pi(T)]\subseteq\pgc(N_\ell\rtimes E)$—as desired.
\end{proof}

We now have sufficient tools to effectively prove \thref{structuralresult} for edges between $\pi(T)$ and $\pi(G)\setminus\pi(T)$ in the prime graph complements of the groups considered in this subsection.
\begin{proposition}
    \thlabel{externalstructure}
    Let $T$ be a non-abelian simple group with Schur multiplier $1$ or $2$, $G$ be of the form $G\cong N.T$ for $N$ solvable, and $E$ be of maximal order among perfect central extensions of $T$ which appear as sections of $G$. Suppose that furthermore $G$ satisfies at least one of the following:
    \begin{enumerate}
        \item $2-p\in\pgc(G)$ for some $p\in\pi(G)\setminus\pi(T)$.
        \item If $t\in\pi(T)$ is Fermat, then for all odd $q \in \pi(G)$ such that $t-q\in\pgc(G)$, there exists a non-Fermat $r\in\pi(T)$ such that $r-q\in\pgc(G)$.
    \end{enumerate}
    Then there exist solvable groups $B$ and $K$ such that $B$ is a direct product of elementary abelian groups, $\pi(K)\cap(\pi(T)\cup N(\pi(T)))=\emptyset$, and $\pi(B)\cup\pi(K)\subseteq \pi(N)$. Furthermore, $E$ and $K$ act by automorphisms on $B$ such that $\pgc(G)\subseteq\pgc(B\rtimes (E\times K))$ and $\pgc(G)$ is formed by removing edges within $\pi(T)$ from $\pgc(B\rtimes (E\times K))$.
\end{proposition}
\begin{proof}
    We remark that $\pgc(G)[\pi(T)]\subseteq \pgc(E)$. If $G\cong L.E$ for some solvable group $L$, then by \thref{2connectsoutward,solvable-rprime}, every prime in $\pi(T)$ which is adjacent to a prime in $\pi(G)\setminus\pi(T)$ does not divide $|L|$. Otherwise, $E\cong2.T$ and recall we still know $G\cong N.T$. Then by \thref{2connectsoutward}, $2-p\not\in\pgc(G)$ for all $p\in\pi(G)\setminus\pi(T)$. Since an order $2$ element is central in $E$, we can furthermore conclude that $2$ is totally isolated within $\pgc(G)$. By \thref{solvable-rprime}, every prime in $\pi(T)$  adjacent to a prime in $\pi(G)\setminus\pi(T)$ does not divide $|N|$. To handle both of these two cases simultaneously, we write $G$ in the form $G\cong L.E_0$ where $L$ is solvable such that every prime in $\pi(T)$ which is adjacent to a prime in $\pi(G)\setminus\pi(T)$ does not divide $|L|$ and $E_0$ is equal to either $E$ if possible, and $T$ otherwise.

    Partition $\pi(T)$ into $\mathcal S = \pi(T)\setminus\pi(L)$ and $\mathcal R=\pi(L)\cap\pi(T)$. For every $p \in \pi(G)\setminus\pi(T)$ such that $r-p\in\pgc(G)$ for some $r\in\pi(T)$, we want to find a corresponding characteristic $p$ representation $\rho$ of $E$ such that for all $s\in\mathcal S$, $s-p\not\in\pgc(G)$ if and only if elements of order $s$ act with fixed points through $\rho$. Fix such a $p$, if it exists.
    
    If $2-q \in \pgc(G)$ for some $q\in\pi(G)\setminus\pi(T)$ then recall $L$ is of odd order. Thus, by our hypotheses on $G$ we may assume that $r$ is non-Fermat or $L$ is of odd order. For each $s\in S$ such that $s-p\not\in\pgc(G)$, let $H^{(s)}$ be the subgroup of $G$ given by \thref{p-maximal}. Then $H^{(s)}$ has a chief series, $1=H^{(s)}_0 \unlhd H^{(s)}_1 \unlhd \ldots \unlhd H^{(s)}_k \unlhd \ldots \unlhd H^{(s)}$, such that $H^{(s)}/H^{(s)}_k \cong E_0$ and if $\ell$ is the maximal integer such that $H^{(s)}_{\ell}/H^{(s)}_{\ell-1}$ is a $p$-group then $s-p\not\in\pgc(H^{(s)}/H^{(s)}_{\ell-1})$. Our argument proceeds with only the quotient group $H^{(s)}/H^{(s)}_{\ell-1}$, so without loss of generality, suppose $H^{(s)}_{\ell-1} = 1$. By Schur-Zassenhaus, $H^{(s)}\cong H^{(s)}_\ell \rtimes (H_k^{(s)}/H^{(s)}_\ell).E_0$. We apply \thref{flavellapplication} to obtain a section $\hat E$ of $G$ which is a perfect central extension of $T$ containing $E_0$ as a section and an action by automorphisms of $\hat E$ on an elementary abelian $p$-group such that some element of order $s$ acts with a fixed point and all elements of order $t\in\pi(T)$ acts Frobeniusly whenever $t-p \in \pgc(G)$. This action gives a characteristic $p$ representation $\rho_s$ of $\hat E$. If $E_0 \cong E$ then by maximality $\hat E \cong E$. Otherwise, when $E_0 \not\cong E$, we recall $2$ is totally isolated within $\pgc(G)$. Thus, if we extend $\rho_s$ to $E$ by \thref{fixedpointextension}, the extension will satisfy the same fixed-point properties as $\rho_s$. So we see that in both cases for $E_0$ we may assume $\rho_s$ is a characteristic $p$ representation of $E$. Applying \thref{fixedpointunions} to the representations $\rho_s$ for each $s\in\mathcal{S}$ such that $s-p\not\in\pgc(G)$ produces a finite $E$-module of characteristic $p$ such that
    \begin{enumerate}
        \item for all $s\in \mathcal S$, if $s-p\not\in\pgc(G)$ then an element of $E$ with order $s$ acts with a fixed point and
        \item for all $t\in\pi(T)$, if $t-p\in\pgc(G)$ then all elements of $E$ with order $t$ act Frobeniusly.
    \end{enumerate}

    Repeat the above process for all $p\in\pi(G)\setminus\pi(T)$ such that $p$ is adjacent to some prime in $\pi(T)$. For all other $p\in\pi(G)\setminus\pi(T)$, associate the trivial representation of $E$. Let $\Xi$ be a graph on vertex set $\pi(G)$ such that $\pgc(E)\subseteq \Xi$, $\Xi\setminus\pi(T) = \pgc(G)\setminus\pi(T)$, and for each $p\in\pi(G)\setminus\pi(T)$ and each $r\in\pi(E)$ we have $r-p\in\Xi$ if and only if all elements of order $r$ in $G$ act Frobeniusly onto the characteristic $p$ module of $E$ the previous process has identified. Since $L$ is solvable and $\Xi\setminus\pi(E) \subseteq \pgc(L)$, we know by \thref{solvableclassification} that $\Xi\setminus\pi(E)$ is triangle free and by \thref{monochromaticneighbors} we have that $\Xi\setminus\pi(E)$ admits the 3-coloring in the hypotheses of \thref{solvablestructure}. Thus, we may apply \thref{solvablestructure} to find the groups $J$ and $K$ such that $\pgc(J\rtimes (E\times K)) = \Xi$. Immediately we have by our construction that
    \begin{enumerate}
        \item $\pi(K)\cap(\pi(T)\cup N(\pi(T))) = \emptyset$ and $\pi(J)\cup \pi(K)\subseteq\pi(N)$.
        \item $\pgc(J \rtimes (E\times K))[\pi(T)] = \pgc(E)$.
        \item $\pgc(J \rtimes (E\times K))\setminus\pi(T) = \pgc(G)\setminus\pi(T)$.
        \item For all $s\in \mathcal S$ and $p\in\pi(G)\setminus\pi(T)$, $s-p\in\pgc(J \rtimes (E\times K))$ if and only if $s-p\in\pgc(G)$.
    \end{enumerate}
    
    By \thref{removingneighbors}, for each $r\in \mathcal R$ we obtain a section $\hat E$ of $G$ which is a perfect central extension of $T$ containing $E_0$ as a section and an action by automorphisms of $\hat E$ on an elementary abelian $r$-group $A_r$ such that $\pgc(G)[\pi(T)]\subseteq \pgc(A_r\rtimes \hat E)$. If $E_0 \cong E$ then by maximality $\hat E \cong E$. Otherwise, when $E_0 \not\cong E$, we recall $2$ is totally isolated within $\pgc(G)$. Then, consider the action by automorphisms of $\hat E$ on $A_r$ as a modular representation of $\hat E$ in characteristic $r$, and extend it to $E$ by \thref{fixedpointextension}. Since $2-r\not\in\pgc(G)$, we see $\pgc(G)[\pi(T)]\subseteq\pgc(A_r\rtimes E)$. So in both cases for $E_0$ we may assume $A_r$ is acted on by $E$.
    
    Let $A= \bigtimes_{r\in\mathcal R} A_r$ and let $E$ act on $A$ as it does on its components. Since no two elements of $\mathcal R$ are adjacent by \thref{divisibilityadjacency}, we have $\pgc(G)[\pi(T)]\subseteq \pgc(A\rtimes E)$. Furthermore, we define the solvable group $B$ by $B=J \times A$ and let $E$ act on $B$ as it does on each of its components. Notice the following all hold:
    \begin{enumerate}
        \item $\pgc(G)[\pi(T)]\subseteq\pgc(B \rtimes (E\times K))[\pi(T)]$.
        \item Since $A$ is an $\mathcal S^\prime$-group, for all $s\in \mathcal S$ and $p\in\pi(G)\setminus\pi(T)$, $s-p\in\pgc(B \rtimes (E\times K))$ if and only if $s-p\in\pgc(G)$.
        \item By \thref{solvable-rprime,2connectsoutward}, for all $r\in \mathcal R$ and $p\in\pi(G)\setminus\pi(T)$, $r-p\not\in\pgc(B \rtimes (E\times K))$ and $r-p\not\in\pgc(G)$.
        \item $\pgc(G)\setminus\pi(T) = \pgc(B \rtimes (E\times K))\setminus\pi(T)$.
    \end{enumerate}
    It thus follows that $\pgc(G)\subseteq\pgc(B \rtimes (E\times K))$, $\pgc(G)$ is formed by removing edges within $\pi(T)$ from $\pgc(B \rtimes (E\times K))$, $\pi(K)\cap(\pi(T)\cup N(\pi(T))) = \emptyset$, and $\pi(K)\cup\pi(B)\subseteq\pi(K)\cup\pi(B_0)\cup\pi(A) \subseteq \pi(N)$.
\end{proof}

To handle edges in the prime graph complement of a group between elements of $\pi(T)$, we will need a version of \thref{flavellapplication} which is able to generate a representation that preserves edge-existence when the normal subgroup is not coprime to both vertices in the edge. The next lemma is designed for this purpose. Similarly to \thref{flavellapplication}, the proof relies on \cite[Theorem A]{Flavell}. Consequently, the assumptions of that theorem impose many of the conditions required by \thref{structuralresult}.
\begin{lemma}
\thlabel{s-action}
    Let $T$ be a non-abelian simple group and $G$ be of the form $G\cong N.E$ with $N$ solvable and $E$ a perfect central extension of $T$. Let $s\in\pi(T)$ and $t\in\pi(T)$ such that $t\nmid|N|$, $t-s\not\in\pgc(G)$, and $t-s\in\pgc(E)$. Suppose that furthermore $G$ satisfies at least one of the following:
    \begin{enumerate}
        \item $s=2$
        \item $N$ is a $2^\prime$-group.
        \item If $u\in\pi(T)$ is Fermat and $u-s\in\pgc(G)$, there exists a non-Fermat $r\in\pi(T)$ such that $r-s\in\pgc(G)$.
    \end{enumerate}
    Then there exists an elementary abelian $s$-group $A$, a section $\hat E$ of $G$ which is a perfect central extension of $T$, and an action by automorphisms of $\hat E$ on $A$ such that $\pgc(G)[\pi(T)]\subseteq\pgc(A\rtimes \hat E)$ and $t-s\not\in\pgc(A\rtimes \hat E)$.
\end{lemma}
The proof of this lemma proceeds in a double-layered induction. The following two lemmas each provide one such inductive step.

\begin{lemma}
\thlabel{subextension}
    Let $T$ be a non-abelian simple group, $s\in\pi(T)$, $E$ a perfect central extension of $T$ by an $s^\prime$-group, $S$ a—possibly trivial—$s$-group, and $G$ a perfect central extension of $S.E$ by an $s^\prime$-group. Furthermore, assume that in no subgroup of $G$ containing $E$ as a section is an order $s$ element central. Then there exists a subgroup $H\leq G$ of the form $H\cong \hat S.\hat E$ where $\hat S$ is an $s$-group and $\hat E$ is a perfect central extension of $T$ containing $E$ as a section.
\end{lemma}
\begin{proof}
    Let $G$ be a counterexample of minimal order. If $Z(G).S\leq \Phi(G)$, then $Z(G).S\cong S\times Z(G)$, as $\Phi(G)$ is nilpotent. Hence $G\cong S.Z(G).E$. Let $S=\hat S$, $\hat E = Z(G).E$, and $H=G$. Since $G$ is perfect, we know $\hat E$ is perfect. Furthermore, by Grün's lemma $\hat E/Z(\hat E)$ has trivial center. Thus $\hat E/Z(\hat E)\cong T$, so $\hat E$ is a perfect central extension of $T$. Then $G$ is not a counterexample.

    Suppose that $Z(G).S\not\leq\Phi(G)$. We know  $Z(G)\leq\Phi(G)$ by \cite[Theorem A.9.3(d)]{frattini}. Then
    \begin{align*}
        \tilde S = (Z(G).S)\Phi(G)/\Phi(G) \cong Z(G).S/(Z(G).S\cap \Phi(G))
    \end{align*}
    is a non-trivial normal $s$-subgroup of $G/\Phi(G)$. By \cite[Theorem A.9.2(e)]{frattini}, $G/\Phi(G)$ has trivial Frattini subgroup, so by \cite[Theorem A.9.2(b)]{frattini} $\tilde S$ has a proper supplement $K/\Phi(G)$. We see $K$ has $E$ as a quotient group. Since $E$ is perfect, if we take $K^{(\infty)}$ to be the perfect core of $K$ (that is, the final term of the derived series), $K^{(\infty)}$ also has $E$ has a quotient group. By assumption, no order $s$-element is central in $K^{(\infty)}$, so $K^{(\infty)}$ matches the hypotheses of the lemma. We recall $K^{(\infty)}$ is a proper subgroup of $G$, so by our minimality hypothesis it has the desired subgroup $H$. Since $H$ is also a subgroup of $G$, we reach our final contradiction.
\end{proof}
\begin{lemma}
\thlabel{s-actioninduction}
    Let $T$ be a non-abelian simple group, $r, s\in\pi(T)$, $G$ be a group of the form $G\cong N.E$ where $N$ is a solvable $r^\prime$-group and $E$ is a perfect central extension of $T$. Suppose that in every section of $G$ containing $T$ as a section no order $r$ or $s$ element is central, and suppose that for every prime $p\in \pi(N) \setminus \{s\}$ an order $r$ element of $G$ acts trivially on some Sylow $p$-subgroup of $N$. Then there exists a subgroup $H\leq G$, of the form $S.\hat E$ where $S$ is an $s$-group and $\hat E$ is a perfect central extension of $T$.
\end{lemma}
\begin{proof}
    Let $G$ be a counterexample of minimal order. If there are no abelian chief factors for $G$ of  characteristic distinct from $s$, then $N$ is an $s$-group, so taking $N=S$ and $E=\hat E$ shows $G$ is not a counterexample. Take a chief series for $G$: $1=N_0\unlhd N_1 \unlhd \ldots \unlhd N_k = N\unlhd \ldots \unlhd E$, and let $i$ be the maximal integer such that $N_i/N_{i-1}$ is an abelian chief factor with characteristic $p$ not equal to $s$. By assumption, let $P$ be the Sylow $p$-subgroup acted trivially upon by an order $r$ element $x \in G$. Then, for some $g\in G$, we have $N_i/N_{i-1}$ is a section of $g^{-1}Pg$. Identifying $g^{-1}xg$ with its corresponding order $r$ element of $G/N_{i-1}$, we see that an order $r$ element of the quotient group centralizes $N_i/N_{i-1}$. Let $C = \mathbf{C}_{G/N_{i-1}}(N_i/N_{i-1})$. We see
    \begin{align*}
        C/((N/N_{i-1})\cap C) \cong \faktor{C(N/N_{i-1})}{(N/N_{i-1})} \unlhd \faktor{(G/N_{i-1})}{(N/N_{i-1})}\cong G/N \cong E.
    \end{align*}
    Recall that $r$ does not divide the order of $N$ but does divide the order of $C$, and that by assumption $Z(E)$ is an $r^\prime$-group. The only normal subgroups of $E$ are central or $E$ itself. Thus, $C/(N/N_{i-1}\cap C) \cong E$ and $C\cong (N_i/N_{i-1}).S_0.E$, where $S_0$ is an $s$-group. Since $N_i/N_{i-1}$ is central, $C/Z(C)$ must be isomorphic to a quotient group of $S_0.E$.
    
    Let $C^{(\infty)}$ denote the perfect core of $C$ (that is, the final term in the derived series). Since $E$ is perfect, $E$ is a quotient of $C^{(\infty)}$. Notice,
    \begin{align*}
        C^{(\infty)}/Z(C^{(\infty)}) \cong \faktor{(C^{(\infty)}/(Z(C)\cap C^{(\infty)}))}{(Z(C^{(\infty)})/(Z(C)\cap C^{(\infty)}))}.
    \end{align*}
    Focusing on the numerator shows
    \begin{align*}
        (C^{(\infty)}/(Z(C)\cap C^{(\infty)})) \cong C^{(\infty)}Z(C)/Z(C) \unlhd C/Z(C).
    \end{align*}
    So $C^{(\infty)}/Z(C^{(\infty)})$ is isomorphic to a normal subgroup of a quotient group of $S_0.E$. Since, by assumption, no order $s$ element of $C^{(\infty)}$ in $Z(C^{(\infty)})$, we then have $C^{(\infty)}$ satisfies the hypotheses of $G$ in \thref{subextension}. Let $K\leq C^{(\infty)}$ be the corresponding subgroup and let $\tilde K$ be the lift of $K$ mod $N_{i-1}$. We see $\tilde K\cong \hat N.\hat E$ where $\hat N$ is a solvable $r^\prime$-group and $\hat E$ is a perfect central extension of $T$. We see furthermore that $\hat N \cong N_{i-1}.S$ where $S$ is an $s$-group. Let $q \in \pi(\hat N)\setminus\{s\}$. Then $q\in\pi(N_{i-1})\setminus\{s\}$. Thus, every Sylow $q$-subgroup of $\hat N$ is an intersection with $N_{i-1}$ of a Sylow $q$-subgroup of $N$. Let $Q$ be a Sylow $q$-subgroup of $N$ acted trivially on by an order $r$ element $y\in G$. Since $r\nmid|N|$ and $r$ is not central in $E$ by assumption, there exists some $g\in G$ so $g^{-1}yg \in \tilde K$. Then $g^{-1}yg$ is an order $r$ element of $\tilde K$ which acts trivially on $g^{-1}Qg \cap N_{i-1}$—a Sylow $q$-subgroup of $\hat N$. So $\tilde K$ satisfies the hypotheses of the lemma. If $\tilde K\neq G$ then $|\tilde K| < |G|$. So, by minimality, it has the desired subgroup $H$. Since $H\leq G$, this is a contradiction. So we can assume $\tilde K=G$.

    Replace $E$ by $\hat E$ and $N$ by $\hat N$ and repeat the above argument. We see that at every iteration, one chief factor is removed from $N$. So, in a finite number of steps, we will have that $G \cong S.\hat E$—a contradiction.
\end{proof}

Now we are able to continue with the proof of \thref{s-action}.
\begin{proof}[Proof of \thref{s-action}]
    Suppose in $\pgc(G)$ that $s$ has no neighbors in $\pi(T)$. Then $A=C_s$ and $E$ acting trivially on $A$ suffices. Now suppose $r-s\in\pgc(G)$ for some $r\in\pi(T)$. By assumption, we may assume $s=2$, $N$ is a $2^\prime$-group, or $r$ is non-Fermat. We know $t\nmid|N|$, there exists an element of $G$ with order $st$, and no such element exists in $E$. Thus, $s\mid|N|$. Furthermore, by \thref{divisibilityadjacency} we then have, $r\nmid|N|$. Take a chief series for $G$ through $N$, $1=N_0\unlhd N_1 \unlhd \ldots \unlhd N_k=N \unlhd \ldots\unlhd G$, and let $\ell$ be the maximal integer such that $s-t\not\in\pgc(G/N_{\ell-1})$. We know $N_\ell/N_{\ell-1}$ is an elementary abelian $s$-group and by \cite[Corollary 1.6]{Suz}, we have $\pgc(G)[\pi(T)]\subseteq\pgc(N_{\ell}/N_{\ell-1}\rtimes G/N_{\ell})$ such that $s-t\not\in\pgc(N_{\ell}/N_{\ell-1}\rtimes G/N_{\ell})$. We will complete the proof within $G/N_{\ell-1}$—so without loss of generality suppose $N_{\ell-1} = 1$. Furthermore, if $K/N_{\ell}$ is the kernel of the action of $G/N_{\ell}$ on $N_{\ell}$ then 
    \begin{align*}
        G/(K\cap N) \cong \faktor{(G/N_{\ell})}{(K/N_{\ell}\cap N/N_{\ell})},
    \end{align*}
    so $\pgc(G)[\pi(T)] \subseteq \pgc(N_{\ell}\rtimes G/(K\cap N))$ such that $s-t\not\in\pgc(N_{\ell}\rtimes (G/(K\cap N))$. We know $G/(K\cap N)$ has $E$ as a quotient group. So again without loss of generality, suppose $K\cap N = 1$ and thus that $N/N_{\ell}$ acts faithfully on $N_\ell$. 
    
    Let $L = N/N_{\ell}$ so $G/N_{\ell} \cong L.E$. Let $p$ be any prime in $\pi(L)$ not equal to $s$ and $R$ be any subgroup of $L.E$ isomorphic to $C_r$. Since $r-s\in\pgc(G)$, $R$ acts faithfully on $N_\ell$. As $(r, |L|) = 1$ there exists an $R$-invariant Sylow $p$-subgroup $P$ of $L$. Consider the corresponding subgroup $P\rtimes R\leq L.E$. Since $L$ and $R$ act faithfully on $N_\ell$ and $R$ acts Frobeniusly on $N_\ell$, we know $R$ acts trivially on $P$ by \cite[Lemma 2.1.6]{2023REU}. Furthermore, since $r-s \in \pgc(G)$, all elements of orders $r$ and $s$ must be non-central in every section of $L.E$ with $T$ as a section. We apply \thref{s-actioninduction} to find a subgroup $H\leq L.E$ of the form $S.\hat E$, where $S$ is an $s$-group and $\hat E$ is a perfect central extension of $T$. Since $t-s\not\in\pgc(G)$ but $t-s\in\pgc(L.E)$, there exists an order $t$ element $x\in G$ which acts with a fixed point on $N_\ell$. Since $t-s\in\pgc(E)$, we know $t\nmid|Z(E)|$, and we recall $t\nmid|N|$. Thus, every Sylow $t$-subgroup of $H$ is a Sylow $t$-subgroup of $L.E$, so there exists $g\in G$ such that $g^{-1}xg \in H$. Then, $t-s\not\in\pgc(N_\ell \rtimes H)$, and the lemma now follows with $A$ as given by \cite[Lemma 3.5]{Suz} and $\hat E$ as defined previously.
\end{proof}

The following proposition completes this subsection by handling the edges within $\pi(T)$. The proof strategy is similar to \thref{externalstructure}, except for some differences in the supporting lemmas cited (see, for example, the discussion above \thref{s-action}). Another notable difference in these two proofs—which may be of note to readers wishing to generalize \thref{structuralresult}—is discussed below the proof of the proposition.
\begin{proposition}
\thlabel{internalstructure}
    Let $T$ be a non-abelian simple group with Schur multiplier $1$ or $2$, $G$ be of the form $G\cong N.T$, and $E$ be of maximal order among perfect central extensions of $T$ which appear as sections of $G$. Suppose that furthermore $G$ satisfies at least one of the following:
    \begin{enumerate}
        \item $2-p\in\pgc(G)$ for some $p\in\pi(G)\setminus\pi(T)$.
        \item If $t\in\pi(T)$ is Fermat, then for all odd $q \in \pi(G)$ such that $t-q\in\pgc(T)$, there exists a non-Fermat $r\in\pi(T)$ such that $r-q\in\pgc(G)$.
    \end{enumerate}
    Then there exists a direct product of elementary abelian groups $B$ and an action by automorphisms of $E$ on $B$ such that $\pgc(G)[\pi(T)] = \pgc(B\rtimes E)$, $\pi(B)\subseteq\pi(N)$, and $r-p\not\in\pgc(G)$ for all $r\in\pi(B)$ and all $p\in\pi(G)\setminus\pi(T)$.
\end{proposition}
\begin{proof}
    We remark that $\pgc(G)[\pi(T)]\subseteq \pgc(E)$. If $G\cong L.E$ for some solvable group $L$, then by \thref{2connectsoutward,solvable-rprime}, every prime in $\pi(T)$ which is adjacent to a prime in $\pi(G)\setminus\pi(T)$ does not divide $|L|$. Otherwise, $E=2.T$ and recall we still know $G\cong N.T$. Then by \thref{2connectsoutward}, $2-p\not\in\pgc(G)$ for all $p\in\pi(G)\setminus\pi(T)$. Thus, by \thref{solvable-rprime}  every prime in $\pi(T)$ which is adjacent to a prime in $\pi(G)\setminus\pi(T)$ does not divide $|N|$. To handle both of these two cases simultaneously, we write $G$ in the form $G\cong L.E_0$, such that $L$ is solvable, every prime in $\pi(T)$ which is adjacent to a prime in $\pi(G)\setminus\pi(T)$ does not divide $|L|$, and $E_0$ is equal to either $E$ if possible, and $T$ otherwise.

    Partition $\pi(T)$ into $\mathcal T = \pi(T)\setminus\pi(L)$ and $\mathcal S=\pi(L)\cap\pi(T)$. For each $s\in\mathcal S$, we wish to find an elementary abelian $s$-group $A_s$, and an action by automorphisms of $E$ on $A_s$ such that $\pgc(G)[\pi(T)]\subseteq\pgc(A_s\rtimes E)$ and $t-s\in\pgc(A_s\rtimes E)$ if and only if $t-s\in\pgc(G)$ for all $t\in\mathcal{T}$. Fix some $s\in\mathcal S$. Recall that if $2-p\in\pgc(G)$ for some $p\in\pi(G)\setminus\pi(T)$ then $L$ is of odd order. So we may always assume that either $L$ is of odd order, $s=2$, or whenever $u-s\in\pgc(G)$ for some Fermat $u \in \pi(T)$ there exists some non Fermat $r\in\pi(T)$ such that $r-s\in\pgc(G)$.
    
    Suppose there exists $s-t\in\pgc(E)$ such that $s-t\not\in\pgc(G)$ and $t\in\mathcal{T}$. We remark that $s-t\in\pgc(E_0)$. By \thref{s-action}, we obtain an elementary abelian $s$-group $A_{s,t}$, a section $\hat E$ of $G$ which is a perfect central extension of $T$, and an action by automorphisms of $\hat E$ on $A_{s,t}$ such that $\pgc(G)[\pi(T)]\subseteq\pgc(A_{s,t}\rtimes \hat E)$ and $s-t\not\in\pgc(A_{s,t}\rtimes \hat E)$. If $\hat E \not\cong E$ then $\hat E \cong T$ and $E\cong2.T$. Then, view the action by automorphisms of $\hat E$ on $A_{s,t}$ as a modular representation of $\hat E$ in characteristic $s$, and extend it to $E$ by \thref{fixedpointextension}. Since $2-s\not\in\pgc(G)$, $\pgc(G)[\pi(T)]\subseteq\pgc(A_{s,t}\rtimes E)$. In both cases for $\hat E$, we may assume $A_{s,t}$ is acted on by $E$. Repeat this process for all $t\in\mathcal T$ such that $t-s\in\pgc(E)$ but $t-s\not\in\pgc(G)$. 
    Let $A_s = \bigtimes_{t\in\mathcal{T}}A_{s,t}$ and $E$ act on $A_s$ as it does on each of its components.

    If no such $t \in\mathcal {T}$ exists, instead apply \thref{removingneighbors} to obtain the elementary abelian group $A_s$, and by an identical argument to above, assume $E$ acts on $A_s$ such that $\pgc(G)[\pi(T)]\subseteq\pgc(A_s\rtimes E)$.

    In either case, we have that for all $t\in\mathcal{T}$, $t-s\in\pgc(A_s\rtimes E)$ if and only if $t-s\in\pgc(G)$ for all $t\in\mathcal{T}$. Repeat the above process for all $s\in\mathcal{S}$, let $B = \bigtimes_{s\in\mathcal{S}}A_s$, and $E$ act on $B$ as it does each of its components. We notice the following:
    \begin{enumerate}
        \item If $u, t \in\mathcal{T}$ then $u-t\in\pgc(G)$ if and only if $u-t\in\pgc(E_0)$. Furthermore, if $E\neq E_0$, then $u,t\neq 2$, so $u-t\in\pgc(G)$ if and only if $u-t\in\pgc(E)$. Finally, since $u,t\nmid|B|$ by construction, $u-t \in \pgc(G)$ if and only if $u-t\in\pgc(B\rtimes E)$.
        \item For $s\in\mathcal S$ and $t\in\mathcal T$, $s-t \in \pgc(G)$ if and only if $s-t\in\pgc(B\rtimes E)$, by construction.
        \item If $s,r\in\mathcal{S}$ then by \thref{divisibilityadjacency}, $s-r\not\in\pgc(G)$ and $s-r\not\in\pgc(A\rtimes E)$.
    \end{enumerate}
    We thus see $\pgc(G)[\pi(T)] = \pgc(B\rtimes E)$. We complete the proof by noticing that $\pi(B)\subseteq\pi(L)$ by construction and thus $\pi(B)\subseteq\pi(N)$ and for all $r\in\pi(B)$ and all $p\in\pi(G)\setminus\pi(T)$, $r-p\not\in\pgc(G)$.
\end{proof}
In the above proof, we need to lift representations of $\hat E$ to those of $E$ using \thref{fixedpointextension} based on possible cases for $\hat E$. A similar process is needed in the proof of \thref{externalstructure}, except in that case, since we know $\hat E$ extends $E_0$, we are able to consider cases for $E_0$ itself—which is necessary to maintain control over the external structure of $\pgc(G)$. These subtle difference might be of note to anyone wishing to extend the results of \thref{structuralresult} to simple groups $T$ of larger Schur multiplier, because it is appears that finer control over how various perfect central extension of $T$ appear within $G$ will be necessary.
\subsection{Induction to Outer Automorphisms}
\label{induction}
This subsection concludes the proof of \thref{structuralresult} by extending \thref{externalstructure,internalstructure} to more arbitrary groups, and then combining them to realize the entire desired prime graph complement.

Extending these lemmas to arbitrary groups requires understanding the effect of outer automorphisms of $T$ have on the edges of the prime graph complement. Since we will always assume that $\operatorname{Out}(T)$ is a $2$-group, its effect can only be the removal of edges from the prime graph complement. Under these assumptions, the following lemma says that the outer automorphisms have no effect on the edges between $\pi(T)$ and $\pi(G)\setminus\pi(T)$.
\begin{lemma}
\thlabel{Onotconnectout}
    Let $T$ be a non-abelian simple group with Schur multiplier $1$ or $2$. Let $G$ be a strictly $T$-solvable group such that there exists some edge within $\pi(T)$. By \cite[Lemma 2.1.1]{2023REU}, $G\cong N.T.O$ where $N$ is solvable and $O\leq \operatorname{Out}(T)$. For all $s\in\pi(T)$ and $p\in\pi(G)\setminus\pi(\Aut(T))$, if $s-p\in\pgc(G)$ then $s-p\in\pgc(N.T)$. 
\end{lemma}
\begin{proof}
    The claim is trivial if $s\nmid|O|$, so suppose otherwise. By \thref{2connectsoutward,solvable-rprime} we write $G\cong L.E.O$ where $L$ is an $s^\prime$-group. Let $S$ be a Sylow $s$-subgroup of $E.O$ and $K\leq G$ be the solvable subgroup given by $L.S$. Identify $S$ with its corresponding isomorphic Sylow $s$-subgroup of $G$. Since $s-p\in\pgc(G)$, by \thref{FrobeniusRef} $H$, the Hall $\{s,p\}$ subgroup of $N.S$, is either Frobenius of type $(s,p)$ or $2$-Frobenius of type $(s,p,s)$. Since $L$ is $s^\prime$, it is easy to see that the latter case is impossible. In the former case, we conclude immediately that $S$ is cyclic or generalized quaternion. Suppose we have $x \in G$ of order $s$ which commutes with $y\in G$ of order $p$. Without loss of generality, let $x\in S$. For some $g\in T.O$, $g^{-1}Sg\cap L.E$ is a Sylow $s$-subgroup of $L.E$. Since we see $y \in L$ and $L \unlhd G$, $g^{-1}yg \in L$ commutes with $g^{-1}xg$. Then, since $S$ is cyclic or generalized quaternion, it contains a unique order $s$-element. So $g^{-1}xg\in L.E$ commutes with $g^{-1}yg \in L.E$—completing the claim.
\end{proof}

The results of Subsection \ref{N.T} involve perfect central extensions in groups $G$ of the form $G\cong N.T$. The next lemma assures that if $G\cong N.T.O$ where $N$ and $O$ are both solvable, then any perfect central extension of $T$ appearing as a section of $G$ is indeed a section of $N.T$—which aligns many results for groups of the form $N.T$ with more general $T$-solvable groups.
\begin{lemma}
\thlabel{ESectionBelow}
    Let $T$ be a non-abelian simple group and $G\cong N.T.O$ where $O$ and $N$ are solvable. Then $E$ is a perfect central extension of $T$ appearing as a section of $G$ if and only if $E$ is a section of $N.T$.
\end{lemma}
\begin{proof}
    If $E$ is a section of $N.T$ then it is immediately a section of $G$. Suppose that $E$ is a section of $G$ and let $H\leq G$ have $M\unlhd H$ such that $H/M \cong E$. By the correspondence theorem, we have $Z\unlhd H$ such that $H/Z\cong T$. Consider
    \begin{align*}
        (H\cap N)Z/Z \unlhd H/Z \cong T.
    \end{align*}
    Since $G$ contains only one composition factor of $T$, $Z$ and $N$ must both be solvable. Thus $(H\cap N)Z\neq H$, so $H\cap N \leq Z$. Since we see
    \begin{align*}
        H/(H\cap N) \cong HN/N \leq N.T/N \cong T \cong H/Z,
    \end{align*}
    it follows that $HN = N.T$ and $H\cap N = Z$. $H\cap N$ is a maximal normal subgroup of $H$. We see $H\cap N \leq H\cap N.T$ so either $H\leq N.T$ or $H\cap N.T = H\cap N$. Since $O$ is solvable and
    \begin{align*}
        H/(H\cap N.T) \cong HN.T/N.T \leq G/N.T \cong O,
    \end{align*}
    the latter case leads to a contradiction. Since $H\leq N.T$, $E$ is a section of $N.T$ by definition.
\end{proof}

When an element of order $2$ is not central in any section of $G$, the outer automorphism group of $T$ can remove edges between elements of $\pi(T)$ from the prime graph complement of $G$. Thus, when \thref{externalstructure,internalstructure} produce actions by automorphisms of $T$, we wish to be able to induce those actions to groups of the form $T.O$ where $O\leq \operatorname{Out}(T)$. The next three lemmas allow for this induction.

The following result is inspired by \cite{aurelstackexchange} and discussion in \cite{jashastackexchange}.
\begin{lemma}
\thlabel{2fixed}
    Let $T$ be a non-abelian simple group and $\rho$ be a representation of $T$ over an arbitrary field with the underlying module $V$. Then every element of order $2$ in $T$ acts with a fixed point through $\rho$.
\end{lemma}
\begin{proof}
    Let $g\in T$ be an element of order $2$. If $V$ has characteristic $2$, we are done by \cite[Proposition 26]{Serre} applied to $\rho \vert_{\langle g\rangle}$, so assume otherwise. Suppose now that $\rho(g)$ has no fixed point. As $\rho(g)^2 =I$, the minimal polynomial of $\rho(g)$ divides $x^2-1 = (x-1)(x+1)$. By assumption, $1\neq -1$, so the minimal polynomial splits as a product of distinct roots. It follows that $\rho(g)$ is diagonalizable and has $\pm 1$ as possible eigenvalues. Since $\rho(g)$ has no fixed points, $1$ is not an eigenvalue. Thus, $\rho(g)$ diagonalizes to $-I$. Then $\rho(g)$ is central in $\rho(T)$, so $g \in Z(T/\ker(\rho))$. Since $T$ is simple, $g\in \ker(\rho)$. Everything in the kernel acts with a fixed point—a contradiction.
\end{proof}

\begin{lemma}
\thlabel{inductionfixedpoints}
    Let $G$ be a group, $N\unlhd G$, and $\theta$ be a representation of $N$ over an arbitrary field $F$. Suppose $\mathcal P \subseteq\pi(N)$ is such that elements of orders from $\mathcal P$ act with fixed points through $\theta$ and all elements of other prime orders act Frobeniusly. Then there exists an $F$-representation $\rho = \operatorname{Ind}_N^G(\theta)$ such that elements of prime orders from $\mathcal{P}$ act with fixed points through $\rho$, and all elements of orders in $\pi(N)\setminus(\pi(G/N)\cup\mathcal P)$ act Frobeniusly.
\end{lemma}
\begin{proof}
    Let $W$ be the underlying module of $\theta$ and define the $FG$-module $V$ as the induced module from $W$ as in \cite[Section 4.3]{Webb}. Define $\rho$ to be the $F$-representation corresponding to the module $V$. Since by \cite[Proposition 4.3.1]{Webb}, $\theta$ is a summand of $\rho\vert_N$, we know that elements of prime orders from $\mathcal P$ act with fixed points through $\rho$. Suppose for the sake of contradiction that an element $x \in G$ of prime order $p \in \pi(N)\setminus(\pi(G/N)\cup\mathcal P)$ acts with a fixed point $v \in V$. Since $p\not\in \pi(G/N)$, $o(xN) \neq p$. Thus, $o(xN) = 1$ and $x\in N$. Again by \cite[Proposition 4.3.1]{Webb}, write $\rho$ as a direct sum of conjugates of $\theta$. Since $x$ acts with a fixed point in $\rho$, it must act with a fixed point in some conjugate $g^{-1}\theta g$ for some $g\in G$. Then, $g^{-1}xg$ acts with a fixed point in $\theta$. However, since $N$ is normal, $g^{-1}xg$ is an order $p\not\in\mathcal P$ element of $N$—which is a contradiction.
\end{proof}

\begin{lemma}
\thlabel{inductionfixedpointscor}
    Let $T$ be a non-abelian simple group with $O\leq\operatorname{Out}(T)$ a—possibly trivial—$2$-group. Then, over any field $F$, for every $F$-representation of $T$, there we have a corresponding $F$-representation of $T.O$ such that there exists an element of prime order $p$ in $T$ acting with fixed points if and only if there exists such an element in $T.O$.
\end{lemma}
\begin{proof}
    Let $\theta$ be a representation of $T$ and $\mathcal{P}\subseteq\pi(T)$ be the prime orders for which there exist elements of $T$ acting with fixed points through $\theta$. By the odd order theorem, $2\mid|T|$, and then by \thref{2fixed}, $2\in\mathcal{P}$. Let $\rho$ be the induced representation of $T$ to $T.O$ given by \thref{inductionfixedpoints}. Then we know for all element orders in $\mathcal{P}$, there exists an element of $T.O$ which acts with a fixed point. Since $\pi(O)\subseteq\{2\}$, no elements of prime orders not in $\mathcal{P}$ acts with fixed points through $\rho$.
\end{proof}

The following two propositions each extend one of \thref{externalstructure,internalstructure} to more general groups using the previous lemmas proven in this subsection. The proof strategies for each are primarily the same, although the former is simpler due to \thref{Onotconnectout}.
\begin{proposition}
    \thlabel{externalstructureinduction}
    Let $T$ be a non-abelian simple group with Schur multiplier $1$ or $2$ and $\operatorname{Out}(T)$ a—possibly trivial—$2$-group. Let $G$ be strictly $T$-solvable with at least one edge within $\pi(T)$ in $\pgc(G)$, and $N\unlhd G$ be as given by \cite[Lemma 2.1.1]{2023REU}. Suppose that $G$ satisfies at least one of the following:
    \begin{enumerate}
        \item $2-p\in\pgc(G)$ for some $p\in\pi(G)\setminus\pi(T)$.
        \item If $t\in\pi(T)$ is Fermat, then for all odd $q \in \pi(G)$ such that $t-q\in\pgc(G)$, there exists a non-Fermat $r\in\pi(T)$ such that $r-q\in\pgc(G)$.
    \end{enumerate}
    Then there exist solvable groups $B$ and $K$ such that $B$ is a direct product of elementary abelian groups, $\pi(K)\cap(\pi(T)\cup N(\pi(T)))=\emptyset$, $\pi(B)\cup\pi(K)\subseteq \pi(N)$, and one of the following is satisfied:
    \begin{enumerate}
        \item $2.T$ is not a section of $G$, $G$ has a quotient group $T.O$ where $O\leq \operatorname{Out}(T)$, and $T.O$ and $K$ act by automorphisms on $B$ such that $\pgc(G)\subseteq\pgc(B\rtimes (T.O\times K))$ and $\pgc(G)$ is formed by removing edges within $\pi(T)$ from $\pgc(B\rtimes (T.O\times K))$.
        \item $2.T$ is a section of $G$ and $2.T$ and $K$ act by automorphisms on $B$ such that $\pgc(G)\subseteq\pgc(B\rtimes (2.T\times K))$ and $\pgc(G)$ is formed by removing edges within $\pi(T)$ from $\pgc(B\rtimes (2.T\times K))$.
    \end{enumerate}
\end{proposition}
\begin{proof}
    Write $G\cong N.T.O$ by \cite[Lemma 2.1.1]{2023REU} where $O\leq \operatorname{Out}(T)$ and $N$ is solvable. By Schreier's Conjecture (known to be true), $O$ is also solvable. Let $E$ be of maximal order among perfect central extensions of $T$ which appear as sections of $G$. By \thref{ESectionBelow}, $E$ must be of maximal order among perfect central extensions of $T$ which appear as sections of $N.T$. We apply \thref{externalstructure} to $N.T$ and let $B_0$ and $K$ be the solvable groups obtained such that $\pgc(N.T)\subseteq\pgc(B_0 \rtimes (E\times K))$ and $\pgc(N.T)$ is formed by removing edges within $\pi(T)$ from $\pgc(B_0 \rtimes (E\times K))$. Since $\pgc(G)\subseteq\pgc(N.T)$ and by \thref{Onotconnectout} $\pgc(G)$ is formed by removing edges within $\pi(T)$ from $\pgc(N.T)$, we have that $\pgc(G)\subseteq\pgc(B_0\rtimes (E\times K))$ and $\pgc(G)$ is formed by removing edges within $\pi(T)$ from $\pgc(B_0\rtimes (E\times K))$. We complete the proof in the following two cases.
    
    {\bf Case 1: $E \cong T$.} Decompose $B_0$ into its Sylow subgroups by $B_0 = \bigtimes B_i$ and consider each as a $T$-module. By \thref{inductionfixedpointscor}, for each $i$ we have an induced module of $T.O$ where elements of prime orders in $T.O$ act with fixed points if and only if there exists an element of the same prime order in $T$ acting with a fixed point on $B_i$. Let $B$ be the direct sum of these modules and $T.O$ act on $B$ as it does each of its components. Then $\pgc(G) \subseteq \pgc(B \rtimes (T.O\times K))$ and $\pgc(G)$ is formed from $\pgc(B\rtimes (T.O\times K))$ by removing edges within $\pi(T)$.

    {\bf Case 2: $E \cong 2.T$.} Letting $B=B_0$, the proposition holds.
\end{proof}

\begin{proposition}
    \thlabel{internalstructureinduction}
    Let $T$ be a non-abelian simple group with Schur multiplier $1$ or $2$ and $\operatorname{Out}(T)$ a—possibly trivial—$2$-group. Let $G$ be strictly $T$-solvable with at least one edge within $\pi(T)$ in $\pgc(G)$, and $N\unlhd G$ be as given by \cite[Lemma 2.1.1]{2023REU}. Suppose that $G$ satisfies at least one of the following:
    \begin{enumerate}
        \item $2-p\in\pgc(G)$ for some $p\in\pi(G)\setminus\pi(T)$.
        \item If $t\in\pi(T)$ is Fermat, then for all odd $q \in \pi(G)$ such that $t-q\in\pgc(G)$, there exists a non-Fermat $r\in\pi(T)$ such that $r-q\in\pgc(G)$.
    \end{enumerate}
    Then there exists a direct product of elementary abelian groups $B$ such that $\pi(B)\subseteq \pi(N)$, $r-p\not\in\pgc(G)$ for all $r\in\pi(B)$ and $p\in\pi(G)\setminus\pi(T)$, and one of the following is satisfied:
    \begin{enumerate}
        \item $2.T$ is not a section of $G$, $G$ has a quotient group $T.O$ where $O\leq \operatorname{Out}(T)$, and $T.O$ acts by automorphisms on $B$ such that $\pgc(G)[\pi(T)]=\pgc(B\rtimes T.O)$.
        \item $2.T$ is a section of $G$ and $2.T$ acts by automorphisms on $B$ such that $\pgc(G)[\pi(T)] = \pgc(B\rtimes 2.T)$.
    \end{enumerate}
\end{proposition}
\begin{proof}
    Write $G\cong N.T.O$ by \cite[Lemma 2.1.1]{2023REU} where $O\leq \operatorname{Out}(T)$ and $N$ is solvable. By Schreier's Conjecture (known to be true), $O$ is also solvable. Let $E$ be of maximal order among perfect central extensions of $T$ which appear as sections of $G$. By \thref{ESectionBelow}, $E$ must be of maximal order among perfect central extensions of $T$ which appear as sections of $N.T$. We apply \thref{internalstructure} to $N.T$ and let $B_0$ be the direct product of elementary abelian groups obtained such that $\pgc(N.T)[\pi(T)] = \pgc(B_0 \rtimes E)$. We complete the proof in the following two cases.
    
    {\bf Case 1: $E \cong T$.} Decompose $B_0$ into its Sylow subgroups by $B_0 = \bigtimes B_i$ and consider each as a $T$-module. By \thref{inductionfixedpointscor} for each $i$ we have an induced module of $T.O$ where elements of prime orders in $T.O$ act with fixed points if and only if there exists an element of the same prime order in $T$ acting with a fixed point on $B_i$. Let $B$ be the direct sum of these modules and $T.O$ act on $B$ as it does on each of its components. Since we see that an edge exists in $\pgc(G)[\pi(T)]$ if and only if it exists in $\pgc(T.O)$ and $\pgc(N.T)[\pi(T)]$, we have $\pgc(G)[\pi(T)] =\pgc(B\rtimes T.O)$.

    {\bf Case 2: $E \cong 2.T$.} The vertex $2$ must be totally isolated within $\pgc(G)$ from other vertices in $\pi(T)$. Thus, since $O$ is a $2$-group, $\pgc(G)[\pi(T)] = \pgc(N.T)[\pi(T)]$. Letting $B=B_0$, the proposition holds.
\end{proof}

All that remains is to combine \thref{externalstructureinduction,internalstructureinduction} to prove \thref{structuralresult}. We restate it here in terms of the notation used in this paper.
\begingroup
\renewcommand{\thetheoremmain}{\ref{structuralresult}}
\begin{theoremmain}
\thlabel{structuralresultrestate}
    Let $T$ be a non-abelian simple group with Schur multiplier $1$ or $2$ and $\operatorname{Out}(T)$ a—possibly trivial—$2$-group. Let $G$ be strictly $T$-solvable with at least one edge within $\pi(T)$ in $\pgc(G)$. Suppose that $G$ satisfies at least one of the following:
    \begin{enumerate}
        \item $2-p\in\pgc(G)$ for some $p\in\pi(G)\setminus\pi(T)$.
        \item If $t\in\pi(T)$ is Fermat, then for all odd $q \in \pi(G)$ such that $t-q\in\pgc(G)$, there exists a non-Fermat $r\in\pi(T)$ such that $r-q\in\pgc(G)$.
    \end{enumerate}
    Then there exist solvable groups $A$ and $K$ such that $A$ is a direct product of elementary abelian groups, $\pi(K)\cap(\pi(T)\cup N(\pi(T))) = \emptyset$, and one the following is satisfied:
    \begin{enumerate}
        \item $2.T$ is not a section of $G$, $G$ has a quotient group $T.O$ where $O\leq \operatorname{Out}(T)$, and $T.O$ and $K$ act by automorphisms on $A$ such that $\pgc(G)=\pgc(A\rtimes(T.O\times K))$.
        \item $G$ has a section $2.T$  and $2.T$ and $K$ act by automorphisms on $A$ such that $\pgc(G) = \pgc(A\rtimes (2.T\times K))$.
    \end{enumerate}
\end{theoremmain}
\endgroup
\begin{proof}
    Let $B_1$ and $K$ be as given by \thref{externalstructureinduction}, $B_2$ be as given by \thref{internalstructureinduction}, and $A=B_1\times B_2$. Let $H$ equal either $T.O$ or $2.T$ depending on which case of \thref{externalstructureinduction,internalstructureinduction} holds for $G$. Furthermore, let $K$ act on $A$ by its given action on $B_1$ and trivially on $B_2$ and $H$ act on $A$ as it does on each of its components. We show  $\pgc(G) = \pgc(A\rtimes (H\times K))$.

    An edge $r-p$ for $r\in\pi(G)$ and $p\in\pi(G)\setminus\pi(T)$ is in $\pgc(B_1\rtimes (H\times K))$ if and only if it is in $\pgc(G)$, by construction. We know $\pgc(A\rtimes (H\times K)) \subseteq\pgc(B_1\rtimes(H\times K))$. Since primes in $\pi(B_2)$ have no neighbors from $\pi(G)\setminus\pi(T)$ in $\pgc(G)$, the trivial actions of $K$ and $B_1$ on $B_2$ do not remove such edges from $\pgc(B_1\rtimes (H\times K))$. Furthermore, since $\pi(B_2)\subseteq\pi(T)$, neither the action of $H$ on $B_2$, nor elements from within $B_2$ remove such edges. Thus, an edge $r-p$ for $r\in\pi(G)$ and $p\in\pi(G)\setminus\pi(T)$ is in $\pgc(A\rtimes (H\times K))$ if and only if it is in $\pgc(G)$.

    An edge $r-s$ for $r,s\in\pi(T)$ is in $\pgc(B_2 \rtimes H)$ if and only if it is in $\pgc(G)$, by construction. We know $\pgc(A\rtimes (H\times K)) [\pi(T)]\subseteq\pgc(B_2\rtimes(H\times K))[\pi(T)]$. By \cite[Lemma 2.1.1]{2023REU}, write $G\cong N.T.O$ and then by \thref{divisibilityadjacency}, if two primes divide $|N|$ we know they are not adjacent in $\pgc(G)$. Thus, since $\pi(B_1) \cup \pi(B_2)\subseteq\pi(N)$, neither the actions of $B_2$ on $B_1$ or elements from within $B_1$ remove edges $r-s$ for $r,s\in\pi(T)$ from $\pgc(B_2\rtimes(H\times K))[\pi(T)]$. Since $\pgc(G)[\pi(T)]\subseteq \pgc(B_1\rtimes(H\times K))[\pi(T)]$, the action of $H$ and $K$ on $B_1$ also does not remove such edges. So an edge $r-s$ for $r,s\in\pi(T)$ is in $\pgc(A \rtimes (H\times K))$ if and only if it is in $\pgc(G)$.
\end{proof}

\section{Applications of the Main Theorem}
\label{classifications}
This section is dedicated to demonstrating the usefulness of \thref{structuralresult}, with a focus towards the problem of classifying the prime graph complements of $T$-solvable groups, for various non-abelian simple groups $T$. We begin with a few generalizations of lemmas in previous work which are useful computationally in many applications of \thref{structuralresult}.

First, generalize \cite[Lemma 6.2]{maslova} from simple groups to arbitrary groups which will be useful computationally in the results to follow.
\begin{lemma}
\thlabel{charactersums}
    Let $G$ be a finite simple group, $F$ be a field of characteristic $p > 0$, $V$ be a absolutely irreducible $GF$-module, and $\beta$ be a Brauer character of $V$. If $g \in G$ is an element of prime order distinct from $p$ then
    \begin{align*}
        \operatorname{dim}\mathbf{C}_V(g) = (\beta_{\langle g \rangle}, 1_{\langle g \rangle}) = \frac{1}{|g|}\sum_{x\in\langle g \rangle}\beta(x).
    \end{align*}
\end{lemma}
\begin{proof}
    Let $\rho$ be the homomorphism corresponding to the module $V$ and let $A=\rho(g) \in GL(n, F)$ where $n=\dim V$. We can view $A$ as instead being in the algebraic closure of $F$, which we denote $K$. Then let $W = K\otimes_{F}V$ be the vector space on which $A$ acts. Since the order of $A$ is $|g|$ which is coprime to $p$, we have that $A$ is diagonalizable  over $K$ by \cite[Corollary 2.1.7]{Webb}.

    Now, we see by \cite[Lemma 3.2.2]{Webb} that $R = \frac{1}{|g|}\sum_{i=1}^{|g|} A^i$ gives the projection operator onto the subspace of $W$ of the fixed points for $A$. Fix a basis such that $A$ is diagonal. Then $A^i$ is diagonal for all $i\in \mathbb{Z}$, so $R$ is diagonal. Since $R$ is a projection map, its eigenvalues must be $0$ or $1$—and the multiplicity of $1$ as an eigenvalue is equal to the dimension of the range of $R$, that is, $\dim\mathbf{C}_W(g)$. Since $\dim\mathbf{C}_W(g)$ is furthermore the multiplicity of $1$ as a root of the characteristic polynomial for $A$—and the characteristic polynomial does not depend on the field extension—we indeed have that the multiplicity of $1$ as an eigenvalue for $R$ equals $\dim\mathbf{C}_V(g) = \dim\mathbf{C}_W(g)$.

    Following the procedure in \cite[Chapter 10]{Webb}, fix a $p$-modular system $(L,T,K)$ such that $L$ is a splitting field for $G$ of characteristic $0$ (which we take to be a subfield of $\mathbb{C}$) with discrete valuation and $T$ the valuation ring of $L$ with a maximal ideal $\langle \pi \rangle$ and $T/\langle \pi\rangle \cong K$. Then let $\ \hat{ }\ $ be an isomorphism from the $p^\prime$ roots of unity in $K$ to the $p^\prime$ roots of unity in $L$ so that for each $x\in G$ we have $\beta(x) = \hat\lambda_1 + \hat\lambda_1 + \ldots + \hat \lambda_{n}$ where $\lambda_1,\ldots, \lambda_{n}$ are the eigenvalues of $\rho(x)$. Take some $1\leq j\leq n$. We see that summing across each entry, we have
    \begin{align*}
        \frac{1}{|g|}\sum_{i=1}^{|g|} (A^i)_{jj} = R_{jj}.
    \end{align*}
    where each $(A^i)_{jj}$ is the $j$th eigenvalue of $A^i$. These sums hold if and only if they hold under $\ \hat{ }\ $. That is, $R_{jj} = 1$ if and only if 
    \begin{align*}
        \frac{1}{|g|}\sum_{i=1}^{|g|} \widehat{(A^i)_{jj}} = 1
    \end{align*}
    and $R_{jj} = 0$ if and only if 
    \begin{align*}
        \frac{1}{|g|}\sum_{i=1}^{|g|} \widehat{(A^i)_{jj}} = 0.
    \end{align*}
    
    Finally, by adding all the values $\widehat{R_{jj}}$ together in characteristic $0$ we have the multiplicity of $1$ as an eigenvalue of $R$. Rearranging this sum gives the desired equation:
    \begin{align*}
        \dim\mathbf{C}_V(g) = \sum_{j=1}^{n} \widehat{R_{jj}} = \sum_{j=1}^{n}\frac{1}{|g|}\sum_{i=1}^{|g|} \widehat{(A^i)_{jj}} = \frac{1}{|g|}\sum_{i=0}^{|g|}\sum_{j=1}^{n}\widehat{(A^i)_{jj}} = \frac{1}{|g|}\sum_{x\in\langle g \rangle}\beta(x).
    \end{align*}
\end{proof}

As a corollary, we get a very useful generalization of \cite[Theorem 3.7]{Suz}, also from simple to arbitrary groups. The original theorem applies the assumption of simplicity only to apply \cite[Lemma 6.2]{maslova}—so our generalization \thref{charactersums} proves the corollary immediately. We also modify the statement of the original theorem to clarify that $Y\subseteq \operatorname{IBr}_p(G)$ should be non-empty and $N$ non-trivial—both of which also follow immediately from the original proof.
\begin{corollary}
    \thlabel{suzgen}
    Let $G$ be a group, and let $p\in \pi(G)$. For each $\chi \in \operatorname{IBr}_p(G)$, let $A_\chi$ be the set of edges $\{p-q\mid \exists g \in G \operatorname{ where} |g| = q, \frac{1}{|g|}\sum_{x\in\langle g \rangle}\chi(g)>0\}$. Given some graph $\Lambda$, we have that $\Lambda$ is realizable as the prime graph complement of a group of the form $N.G$ where $N$ is a non-trivial $p$-group if and only if there is some nonempty subset $Y\subseteq \operatorname{IBr}_p(G)$ such that $\Lambda = \pgc(G) \setminus \bigcup_{\chi \in Y}A_\chi$.
\end{corollary}

\subsection{Applicable Groups}
\label{applicablegroups}
The primary purpose of this subsection is the proof of the following, which demonstrates the applicability of \thref{structuralresult} to a number of large classes of groups.
\begin{theorem}
\thlabel{thmapplicablegroups}
    Let $T$ be equal to one of $A_{11}$, $A_{13}$, $A_{14}$, $A_{15}$, $A_{16}$, $\PSL(2,13)$, $S_6(3)$, $HS$, or $G_2(4)$. Then if $G$ is a $T$-solvable group, either $\pgc(G)$ is triangle free and $3$-colorable or $G$ satisfies the hypotheses of \thref{structuralresult}.
\end{theorem}

We will prove the above theorem in the propositions which follow; however we save the proof for the $T=\PSL(2,13)$ case for Subsection \ref{PSL}, where the prime graph complements of $\PSL(2,13)$-solvable groups are fully classified. The following lemma will mostly prove 
\thref{thmapplicablegroups} for $T\neq \PSL(2,13)\text{ or }G_2(4)$. We will see in Subsection \ref{PSL} that the prime graph complements of $\PSL(2,13)$-solvable groups may have an arbitrary number of triangles, and the representation information for $G_2(4)$ in Appendix \hyperref[info]{A} indicates the same may be true for $G_2(4)$-solvable groups. It is thus easier to handle them separately—although the overall proof strategy remains similar to the rest of the groups considered in this subsection.
\begin{lemma}
    \thlabel{applicabilitycriterion}
    Let $T$ be a non-abelian simple group with Schur multiplier $1$ or $2$ and $\operatorname{Out}(T)$ a—possibly trivial—$2$-group. Suppose furthermore that
    \begin{enumerate}
        \item $\pgc(T)$ is $3$-colorable.
        \item There are two primes $p$ and $q$ neither of which is Fermat, such that every triangle in $\pgc(T)$ contains the edge $p-q$.
        \item If $r\in\pi(T)$ is Fermat and $r-s\in\pgc(T)$ then $s\in\{2,p,q\}$.
        \item For any $T$-solvable group $G$, there exists at most one odd prime which is adjacent to $\pi(G)\setminus\pi(T)$, and it is not a Fermat prime.
    \end{enumerate}
    Then if a group $G$ is $T$-solvable, either $\pgc(G)$ is triangle-free and $3$-colorable or $G$ satisfies the hypotheses of \thref{structuralresult}.
\end{lemma}
\begin{proof}
    Let $G$ be $T$-solvable. If $2-r\in\pgc(G)$ for some $r\in\pi(G)\setminus\pi(T)$, then $G$ satisfies the hypotheses of \thref{structuralresult} immediately. So suppose instead that $2$ is not adjacent to $\pi(G)\setminus\pi(T)$. By \thref{solvableclassification} we may furthermore assume $G$ is strictly $T$-solvable. 
    
    We now wish to show that $\pgc(G)$ is $3$-colorable and every triangle in $\pgc(G)$ contains the edge $p-q$. Take $K\leq G$ as given by \cite[Lemma 2.1.2]{2023REU}. Since $\pgc(G)\subseteq\pgc(K)$, it suffices to prove the claim in $K$, so assume without loss of generality that $G=K$. As $\pgc(T)$ is $3$-colorable and only one prime in $\pi(T)$ connects to $\pi(G)\setminus\pi(T)$, it is immediate from \thref{monochromaticneighbors} that $\pgc(G)$ is $3$-colorable. Furthermore, since only one prime in $\pi(T)$ connects to $\pi(G)\setminus\pi(T)$, we know by \cite[Lemma 5.7]{Suz} that every triangle in $\pgc(G)$ is contained within $\pgc(G)[\pi(T)]\subseteq\pgc(T)$. Thus, since every triangle in $\pgc(T)$ contains the edge $p-q$, the same holds for every triangle in $\pgc(G)$.

    Thus, if $\pgc(G)$ lacks the edge $p-q$, then it is triangle-free and $3$-colorable. Now, we recall by assumption the Fermat primes in $\pi(T)$ connect in $\pgc(G)$ only to $2$, $p$ or $q$. Hence, if $\pgc(G)$ has the edge $p-q$, then $G$ satisfies the hypotheses of \thref{structuralresult}. Thus altogether either $\pgc(G)$ is triangle-free and $3$-colorable or $G$ satisfies the hypotheses of \thref{structuralresult}.
\end{proof}

\begin{proposition}
    Let $G$ be an $A_{11}$-solvable group. Then either $\pgc(G)$ is triangle free and $3$-colorable or $G$ satisfies the hypotheses of \thref{structuralresult}.
\end{proposition}
\begin{proof}
    We verify that $A_{11}$ satisfies the hypotheses of \thref{applicabilitycriterion}. We see in Appendix \hyperref[info]{A} that $A_{11}$ has Schur multiplier $2$, $|\operatorname{Out}(A_{11})|=2$, $\pgc(A_{11})$ is $3$-colorable, every triangle contains the edge $7-11$, and the Fermat primes $3$ and $5$ in $\pi(A_{11})$ are adjacent only to $7$ and $11$. Furthermore, the representation information in Appendix \hyperref[info]{A} shows $3$, $5$, and $7$ all satisfy \cite[Corollary 2.2.6]{2023REU}. Hence, only $11\in\pi(A_{11})\setminus\{2\}$ may be adjacent to any $p\in \pi(G)\setminus\pi(A_{11})$. Thus the proposition holds by \thref{applicabilitycriterion}.
\end{proof}

\begin{proposition}
    Let $n\in\{13, 14, 15, 16\}$ and $G$ be $A_n$-solvable. Then either $\pgc(G)$ is triangle free and $3$-colorable or $G$ satisfies the hypotheses of \thref{structuralresult}.
\end{proposition}
\begin{proof}
    We verify that $A_{n}$ satisfies the hypotheses of \thref{applicabilitycriterion}. We see in Appendix \hyperref[info]{A} that $A_{n}$ has Schur multiplier $2$, $|\operatorname{Out}(A_n)|=2$, $\pgc(A_n)$ is $3$-colorable, every triangle contains the edge $11-13$, and the Fermat primes $3$ and $5$ in $\pi(A_{n})$ are adjacent only to $11$ and $13$.
    
    Suppose for the sake of contradiction that there exists an edge $r-p\in\pgc(G)$ for some $p\in\pi(G)\setminus\pi(A_n)$ and $r\in\pi(A_n)\setminus \{ 2, 13\}$. Then $r-p\in\pgc(K)$ where $K\leq G$ is as given by \cite[Lemma 2.1.2]{2023REU}. Since $A_{13}\leq A_n$ and $\pi(A_{13}) = \pi(A_n)$, there exists $L\leq K$ such that $L\cong N.A_{13}$ for $N$ solvable and $r-p\in \pgc(L)$. By the representation information in Appendix \hyperref[info]{A} for $A_{13}$, we have a contradiction with \cite[Corollary 2.2.6]{2023REU}. Hence, only $13\in\pi(A_{13})\setminus\{2\}$ may be adjacent to any $p\in \pi(G)\setminus\pi(A_{n})$. Thus the proposition holds by \thref{applicabilitycriterion}.
\end{proof}

\begin{proposition}
    Let $G$ be an $S_6(3)$-solvable group. Then either $\pgc(G)$ is triangle free and $3$-colorable or $G$ satisfies the hypotheses of \thref{structuralresult}.
\end{proposition}
\begin{proof}
    We verify that $S_{6}(3)$ satisfies the hypotheses of \thref{applicabilitycriterion}. We see in Appendix \hyperref[info]{A} that $S_{6}(3)$ has Schur multiplier $2$, $|\operatorname{Out}(S_6(3))|=2$, $\pgc(S_6(3))$ is $3$-colorable, every triangle contains the edge $7-13$, and the Fermat primes $3$ and $5$ in $\pi(S_6(3))$ are adjacent only to $7$ and $13$. Furthermore, the representation information in Appendix \hyperref[info]{A} shows $3$, $5$, $7$, and $13$ all satisfy \cite[Corollary 2.2.6]{2023REU}. Hence, no odd primes in $\pi(S_{6}(3))$ may be adjacent to any $p\in \pi(G)\setminus\pi(S_6(3))$. Thus the proposition holds by \thref{applicabilitycriterion}.
\end{proof}

\begin{proposition}
    Let $G$ be an $HS$-solvable group. Then either $\pgc(G)$ is triangle free and $3$-colorable or $G$ satisfies the hypotheses of \thref{structuralresult}.
\end{proposition}
\begin{proof}
    We verify that $HS$ satisfies the hypotheses of \thref{applicabilitycriterion}. We see in Appendix \hyperref[info]{A} that $HS$ has Schur multiplier $2$, $|\operatorname{Out}(HS)|=2$, $\pgc(HS)$ is $3$-colorable, every triangle contains the edge $7-11$, and the Fermat primes $3$ and $5$ in $\pi(HS)$ are adjacent only to $7$ and $11$. Furthermore, the representation information in Appendix \hyperref[info]{A} shows $3$, $5$, $7$, and $11$ all satisfy \cite[Corollary 2.2.6]{2023REU}.  Hence, no odd primes in $\pi(HS)$ may be adjacent to any $p\in \pi(G)\setminus\pi(HS)$. Thus the proposition holds by \thref{applicabilitycriterion}.
\end{proof}

\begin{proposition}
\thlabel{G24applicability}
    Let $G$ be a $G_2(4)$-solvable group. Then either $\pgc(G)$ is triangle free and $3$-colorable or $G$ satisfies the hypotheses of \thref{structuralresult}.
\end{proposition}
\begin{proof}
    We see in Appendix \hyperref[info]{A} that $G_2(4)$ has Schur multiplier $2$, $|\operatorname{Out}(G_2(4))|=2$, $\pgc(G_2(4))$ is $3$-colorable, every triangle contains the edge $7-13$, and the Fermat primes $3$ and $5$ in $\pi(G_2(4))$ are adjacent only to $7$ and $13$. Furthermore, the representation information in Appendix \hyperref[info]{A} shows $3$ and $5$ satisfy \cite[Corollary 2.2.6]{2023REU} so the only primes in $\pgc(G)$ which may be adjacent to a Fermat prime are $7$ and $13$. Proceeding similarly to \thref{applicabilitycriterion}, it suffices to prove $\pgc(G)$ is $3$-colorable and every triangle contains the edge $7-13$.

    By \thref{solvableclassification} we may assume $G$ is strictly $G_2(4)$-solvable. Then, if we take $K\leq G$ as given by \cite[Lemma 2.1.2]{2023REU} since $\pgc(G)\subseteq\pgc(K)$ it is sufficient to prove the claim in $K$. So without loss of generality we may assume $G=K$. Furthermore, if $2-r\in\pgc(G)$ for $r\in\pi(G)\setminus\pi(G_2(4))$ then $G$ satisfies the hypotheses of \thref{structuralresult} immediately. So suppose instead that $2$ is not adjacent to $\pi(G)\setminus\pi(G_2(4))$. Color $\pgc(G)\setminus\pi(G_2(4))$ by \thref{monochromaticneighbors}. If we say $N(\pi(G_2(4)))\setminus\pi(G_2(4))$ is colored $\mathcal{I}$ then color $2$, $3$, and $5$ by $\mathcal I$, $7$ by $\mathcal{O}$ and $13$ by $\mathcal{D}$ to give a proper $3$-coloring of $\pgc(G)$. 
    
    Now, take a triangle in $\pgc(G)$. By \cite[Lemma 5.7]{Suz}, at most one vertex of the triangle may be in $\pi(G)\setminus\pi(G_2(4))$. If the triangle is not contained within $\pi(G_2(4))$ then it must contain the edge $7-13$ because only $7$ and $13$ in $\pi(G_2(4))$ can connect to $\pi(G)\setminus\pi(G_2(4))$. Otherwise, it must contain the edge $7-13$ by our earlier observation on the triangles in $\pgc(G_2(4))$.
\end{proof}

One notices that we have shown that all $A_n$-solvable groups with $n\in\{11,$ $ 13,$ $ 14,$ $ 15, $ $16\}$ either have triangle free and $3$-colorable prime graph complements or satisfy the hypotheses of \thref{structuralresult}. We do not consider $A_n$ for larger values of $n$ because we note that $\pgc(A_{17})$ itself is not triangle-free and fails to satisfy the hypotheses of \thref{structuralresult} because $3$ is adjacent only to the Fermat prime $17$. For $5\leq n \leq 10$, the prime graph complements of $A_n$-solvable groups have already been classified in \cite{2021REU,2022REU,2023REU}. We have also excluded $A_{12}$-solvable groups because it is simpler to instead prove the following classification.

\begin{proposition}
    A graph $\Xi$ is isomorphic to the prime graph complement of an $A_{12}$-solvable group if and only if it is triangle free and $3$-colorable.
\end{proposition}
\begin{proof}
    Any triangle free and $3$-colorable graph $\Xi$ is realizable by an $A_{12}$-solvable group, because by \thref{solvableclassification} it is realizable by a solvable group. Now, suppose some group $G$ is $A_{12}$-solvable. If it is not strictly solvable, then the proposition follows by \thref{solvableclassification}, so assume $G$ is strictly $A_{12}$-solvable. We may without loss of generality replace $G$ by the subgroup $K\leq G$ given by \cite[Lemma 2.1.2]{2023REU} since we recall $\pgc(G)\subseteq\pgc(K)$. We see by inspection in Appendix \hyperref[info]{A} that $\pgc(G)[\pi(A_{12})]\subseteq\pgc(A_{12})$ is triangle free and $3$-colorable. By computing the structure description of the Sylow $2$-subgroup of $A_{12}$ in \cite{GAP}, we can verify that $2$ satisfies the hypotheses of \cite[Proposition 2.2.2]{2023REU}, so $2-p\not\in\pgc(G)$ for all $p\in\pi(G)\setminus\pi(A_{12})$. Similarly, the representation information in \hyperref[info]{A} shows $3$, $5$, $7$, and $11$ all satisfy \cite[Corollary 2.2.6]{2023REU}. Hence, no prime in $\pi(A_{12})$ may be adjacent to any $p\in \pi(G)\setminus\pi(A_{12})$. Then $\pgc(G)$ is triangle free by \cite[Lemma 5.7]{Suz}, and $\pgc(G)$ is $3$-colorable by \thref{monochromaticneighbors}.
\end{proof}
\subsection{$\PSL(2,13)$}
\label{PSL}
The following lemma follows quickly from existing edge non-existence criteria.
\begin{lemma}
\thlabel{PSLedgenonexistence}
    Let $G$ be a strictly $\PSL(2,13)$-solvable group. Then $3-p\not\in\pgc(G)$ for all $p\in\pi(G)\setminus\pi(\PSL(2,13))$. Furthermore, if $2-r\in\pgc(G)$ for any $r\in\pi(\PSL(2,13))$ then $2-p\not\in\pgc(G)$ for all $p\in\pi(G)\setminus\pi(\PSL(2,13))$.
\end{lemma}
\begin{proof}
    We see in Appendix \hyperref[info]{A} that that the prime $3$ satisfies \cite[Corollary 2.2.6]{2023REU}. Thus $3-p\not\in\pgc(G)$ for all $p\in\pi(G)\setminus\pi(\PSL(2,13))$.
    
    Now, suppose that $2-p$ is in $\pgc(G)$ for some $p\in\pi(G)\setminus\pi(T)$. Then $2-p \in \pgc(K)$ where $K\leq G$ is given in \cite[Lemma 2.1.2]{2023REU}. By \thref{2connectsoutward}, $2.\PSL(2,13)$ appears as a section of $K$, so $2-r\not\in\pgc(K)$ for all $r\in\pi(\PSL(2,13))$. Hence if $2-r\in\pgc(G)$, then no such edge $2-p$ may exist.
\end{proof}

Now we can use \thref{monochromaticneighbors} to show that $\PSL(2,13)$-solvable groups admit $3$-colorings of a similar form.
\begin{lemma}
    \thlabel{PSLcoloring}
    Let $G$ be a $\PSL(2,13)$-solvable group. Then $\pgc(G)$ is $3$-colorable. Furthermore, if $G$ is strictly $\PSL(2,13)$-solvable, then our coloring may be chosen such that $N(\pi(\PSL(2,13))\setminus\pi(\PSL(2,13))$ is monochromatic.
\end{lemma}
\begin{proof}
    If $G$ is solvable, then the lemma is handled by \thref{solvableclassification}, so suppose $G$ is strictly $\PSL(2,13)$-solvable. It is enough to put the desired coloring on the subgroup $K\leq G$ as given by \cite[Lemma 2.1.2]{2023REU}. Without loss of generality, assume $G=K$. Then we apply \thref{monochromaticneighbors} to color $\pgc(G)\setminus\pi(\PSL(2,13))$. If we say that $N(\pi(\PSL(2,13)))\setminus\pi(\PSL(2,13))$ is colored $\mathcal I$, then color $13$ by $\mathcal O$, $7$ by $\mathcal{D}$, $3$ by $\mathcal{I}$ and $2$ by either $\mathcal{D}$ or $\mathcal{I}$ depending on whether or not it is adjacent to any $p\in\pi(G)\setminus\pi(\PSL(2,13))$, respectively. By inspection of $\pgc(\PSL(2,13))$ in Appendix \hyperref[info]{A} and the restrictions placed in \thref{PSLedgenonexistence}, the described process gives the desired proper $3$-coloring on $\pgc(G)$.
\end{proof}

\begin{lemma}
\thlabel{PSLtriangleedge}
    Let $G$ be a $\PSL(2,13)$-solvable group. Then every triangle in $\pgc(G)$ contains the edge $7-13$.
\end{lemma}
\begin{proof}
    Suppose $\pgc(G)$ contains a triangle. We see that $G$ must be strictly $\PSL(2,13)$-solvable by \thref{solvableclassification}. By \cite[Lemma 5.7]{Suz} applied to the subgroup of $G$ given by \cite[Lemma 2.1.2]{2023REU}, at most one vertex of the triangle may be from $\pi(G)\setminus\pi(\PSL(2,13))$. If the triangle is not contained in $\pi(\PSL(2,13)$ then \thref{PSLedgenonexistence} guarantees it contains the $13-7$ edge. Otherwise, it must be one of the two triangles in $\pgc(\PSL(2,13))$, both of which we can see contain the edge $7-13$ by inspecting Appendix \hyperref[info]{A}.
\end{proof}

We are now able to complete the proof of \thref{thmapplicablegroups}. 
\begin{proposition}
\thlabel{PSLsatisfies}
    Let $G$ be a $\PSL(2,13)$-solvable group. Then either $\pgc(G)$ is triangle free and $3$-colorable, or $G$ satisfies the hypotheses of \thref{structuralresult}. 
\end{proposition}
\begin{proof}
    We see in Appendix \hyperref[info]{A} that $\PSL(2,13)$ has Schur multiplier $2$, and $|\operatorname{Out(\operatorname{PSL}(2,13))}|=2$. Furthermore, by \thref{PSLcoloring}, $\pgc(G)$ is $3$-colorable. So it only remains to consider if $\pgc(G)$ has a triangle. In this case, by \thref{solvableclassification}, we know $G$ is strictly $\PSL(2,13)$-solvable.
    
    Since $3$ is never adjacent to $\pi(G)\setminus\pi(\PSL(2,13))$ by \thref{PSLedgenonexistence}, we see from inspecting $\pgc(\PSL(2,13))$ in Appendix \hyperref[info]{A} that the only possible neighbors of $3$ are the primes $7$ and $13$. Since the edge $7-13$ is present, $3$ is not the sole neighbor of any vertex. Noting that $3$ is the only Fermat prime in $\pi(\PSL(2,13))$, $G$ satisfies the assumptions of \thref{structuralresult}.
\end{proof}

We see below how the strict structure imposed by \thref{structuralresult} makes simple a number of lemmas restricting the structure of the prime graph complements of $\PSL(2,13)$-solvable groups.
\begin{lemma}
    \thlabel{PSL2noedge7noedge}
    Let $G$ be a $\PSL(2,13)$-solvable group satisfying the hypotheses of \thref{structuralresult}. If $2-r \in \pgc(G)$ for some $r \in  \pi(\PSL(2,13))$, then $7-p \not\in\pgc(G)$ for all $p \in \pi(G)\setminus\pi(\PSL(2,13))$.
\end{lemma}
\begin{proof}
    Without loss of generality we may assume $G$ is of the form given by \thref{structuralresult}. Since $2-r\in\pgc(G)$ for some $r\in\pi(\PSL(2,13))$, we see that $2.\PSL(2,13)$ is not a section of $G$, so we assume specifically that $G=A\rtimes (H\times K)$ where $A$, $H$, and $K$ are as given by \thref{structuralresult}.

    Let $p\in\pi(G)\setminus\pi(\PSL(2,13)$ be arbitrary and suppose for the sake of contradiction that $7-p\in\pgc(G)$. Let $K\leq G$ be the subgroup of the form $P\rtimes \PSL(2,13)$ where $P$ is the Sylow $p$-subgroup $A$ and the action of $\PSL(2,13)$ on $P$ is that of the isomorphic subgroup of $H$. 
    
    We see that $7-p\in\pgc(K)$, so every element of order $7$ in $\PSL(2,13)$ must act Frobeniusly on $P$. According to \cite[Lemma 2.1.7]{2023REU}, we then have a $\PSL(2,13)\mathbb{C}$-module $\tilde P$ on which every element of order $7$ acts Frobeniusly. By reviewing the representation information of $\PSL(2,13)$ in Appendix \hyperref[info]{A} we see that no such representation exists—forcing a contradiction.
\end{proof}

\begin{lemma}
\thlabel{PSL7subset}
    Let $G$ be a $\PSL(2,13)$-solvable group $G$ satisfying the hypotheses of \thref{structuralresult}. Suppose $7-p \in\pgc(G)$ for some $p \in \pi(G)\setminus \pi(\PSL(2,13))$. If $13-q\in\pgc(G)$ for some $q\in\pi(G)\setminus\pi(\PSL(2,13))$, then $13-p \in \pgc(G)$. If $2-q\in\pgc(G)$ for some $q\in\pi(G)\setminus\pi(\PSL(2,13))$, then $13-p \in \pgc(G)$.
\end{lemma}
\begin{proof}
    The proofs for $13$ and $2$ are identical, so we will suppose only that $13-q\in\pgc(G)$ for some $q\in\pi(G)\setminus\pi(\PSL(2,13))$.

    Without loss of generality we may assume $G$ is of the form given by \thref{structuralresult}. That is, either $G=A\rtimes(H\times K)$ or $G=A\rtimes (E\times K)$ where $A, K, H,$ and $E$ are as given by \thref{structuralresult}. We see that since $A$ is abelian if $13\mid|A|$ then $13-q\not\in\pgc(G)$. So if $P$ is the Sylow $p$-subgroup of $A$ then it is enough to show $13-p\in\pgc(P\rtimes H)$ or $13-p\in\pgc(P\rtimes E)$—depending on which case for $G$ holds. Let $\tilde P$ be the corresponding $H\mathbb{C}$- or $E\mathbb{C}$-module given by \cite[Lemma 2.1.7]{2023REU}. We know all elements of order $7$ in $\PSL(2,13)$ act Frobeniusly on $\tilde P$, and it is enough to show the same for all elements of order $13$.
    
    Since the outer automorphism group of $\PSL(2,13)$ is $C_2$ we know that $H$ and $E$ can only be one of $\PSL(2,13)$, $2.\PSL(2,13)$, or $\Aut(\PSL(2,13))$. Now, by reviewing the representation information of these groups in Appendix \hyperref[info]{A}, we see that wherever all elements of order $7$ act Frobeniusly through a complex representation of one of these groups, all elements of order $13$ act Frobeniusly. So since all elements of order $7$ act Frobeniusly on $\tilde P$, all elements of order $13$ must also act Frobeniusly on $\tilde P$. We have shown then that $13-p\in\pgc(G)$.
\end{proof}

\begin{lemma}
\thlabel{PSLinternalstructure}
    Let $G$ be a $\PSL(2,13)$-solvable group satisfying the hypotheses of \thref{structuralresult} such that $\pgc(G)$ contains a triangle. Then $\pgc(G)[\pi(\PSL(2,13))]$ is not equal to any of the following:\vspace{10pt}
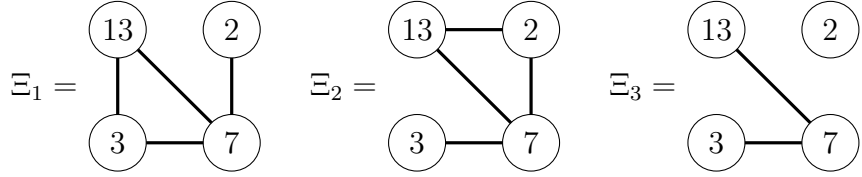
    \begin{center}   
            $\Xi_1=$ \begin{tikzpicture}[baseline=18pt, state/.style={circle, draw, minimum size=0.75cm, inner sep = 0cm, outer sep = 0cm}]
         \node[state] (D) at (0,0) {$3$};
            \node[state] (B) at (1.5,0) {$7$};
            \node[state] (C) at (0,1.5) {$13$};
            \node[state] (A) at (1.5,1.5) {$2$};
            
            \draw [line width=1.2pt](A) -- (B);
            \draw [line width=1.2pt](B) -- (D);
            \draw [line width=1.2pt](D) -- (C);
            \draw [line width=1.2pt](B) -- (C);
        
            \end{tikzpicture}
                \hspace{10pt}
            $\Xi_2=$ \begin{tikzpicture}[baseline=18pt, state/.style={circle, draw, minimum size=0.75cm, inner sep = 0cm, outer sep = 0cm}]
         \node[state] (D) at (0,0) {$3$};
            \node[state] (B) at (1.5,0) {$7$};
            \node[state] (C) at (0,1.5) {$13$};
            \node[state] (A) at (1.5,1.5) {$2$};
            
            \draw [line width=1.2pt](A) -- (B);
            \draw [line width=1.2pt](B) -- (D);
            \draw [line width=1.2pt](A) -- (C);
            \draw [line width=1.2pt](B) -- (C);
        \end{tikzpicture}
        \hspace{10pt}
        $\Xi_3=$ \begin{tikzpicture}[baseline=18pt, state/.style={circle, draw, minimum size=0.75cm, inner sep = 0cm, outer sep = 0cm}]
         \node[state] (D) at (0,0) {$3$};
            \node[state] (B) at (1.5,0) {$7$};
            \node[state] (C) at (0,1.5) {$13$};
            \node[state] (A) at (1.5,1.5) {$2$};
            
            \draw [line width=1.2pt](B) -- (D);
            \draw [line width=1.2pt](B) -- (C);
        \end{tikzpicture}
        \end{center}
\end{lemma}
\begin{proof}
    Without loss of generality we may assume $G$ is of the form given by \thref{structuralresult}. That is, either $G=A\rtimes(H\times K)$ or $G=A\rtimes (E\times K)$ where $A, K, H,$ and $E$ are as given by \thref{structuralresult}.
    
    {\bf Case 1:} Assume for the sake of contradiction that $\pgc(G)[\pi(\PSL(2,13))] = \Xi_1$. Since $2-7\in\pgc(G)$ we see that $2.\PSL(2,13)$ is not a section of $G$. Thus we may assume specifically that $G=A\rtimes(H\times K)$. Furthermore since $2-7\not\in\pgc(\Aut(\PSL(2,13))$, we may assume $H=\PSL(2,13)$. Suppose towards a contradiction that $13\mid|A|$. Then let $W$ be the Sylow $13$-subgroup of $A$ and notice that $7-13\in\pgc(A\rtimes \PSL(2,13))$. By reviewing the mod $13$ Brauer characters of $\PSL(2,13)$ in Appendix \hyperref[info]{A} we have a contradiction from \thref{suzgen}. So we may assume $13\nmid|A|$. Thus if $V$ is the Sylow $2$-subgroup of $A$ we must have $\pgc(V\rtimes \PSL(2,13)) = \Xi_1$. By reviewing the mod $2$ Brauer characters of $\PSL(2,13)$ in Appendix \hyperref[info]{A} we again see a contradiction by \thref{suzgen}.
    
    {\bf Case 2:} This case follows identically to the former except we must review the mod $3$ Brauer characters of $\PSL(2,13)$ rather than the mod $2$ ones.

    {\bf Case 3:} Assume for the sake of contradiction that $\pgc(G)[\pi(\PSL(2,13))] = \Xi_3$. Since we assume $\pgc(G)$ contains a triangle, we know by \thref{PSLtriangleedge} that $7-q\in\pgc(G)$ and $13-q\in\pgc(G)$ for some $q\in\pi(G)\setminus\pi(\PSL(2,13))$. Thus since $A$ is abelian we have that $7\nmid|A|$ and $13\nmid|A|$. Suppose for the sake of contradiction that $G=A\rtimes (H\times K)$. Then if $P$ is the Sylow $p$-subgroup of $A$ we see $7-p\in\pgc(P\rtimes H)$. Thus, every order $7$ element of $H$ must act Frobeniusly on $P$, so on the corresponding $H\mathbb{C}$-module $\tilde P$ from \cite[Lemma 2.1.7]{2023REU}. By reviewing the representation information of $\PSL(2,13)$ and $\Aut(\PSL(2,13))$—the two possibilities for $H$—in Appendix \hyperref[info]{A} we see a contradiction since no such module exists. Thus we may assume $G = A\rtimes(E\times K)$. Now, let $V$ be the Sylow $3$-subgroup of $A$. Recalling that $13\nmid|A|$, we must have that $\pgc(V\rtimes E) = \Xi_3$. By reviewing the mod $3$ Brauer characters of $E=2.\PSL(2,13)$ in Appendix \hyperref[info]{A} we see a contradiction by \thref{suzgen}. 
\end{proof}

The following full classification for the prime graph complements of $\PSL(2,13)$-solvable groups—while long—follows simply from the previous lemmas. The bulk of the proof is using the restrictions we have already imposed on the prime graph complement to show that it lands into one of the cases in the theorem.
\begin{theorem}\thlabel{pslclassification}
    A graph $\Xi$ is isomorphic to the prime graph complement of some $\PSL(2,13)-$solvable group if and only if one of these conditions holds:
    \begin{enumerate}
        \item $\Xi$ is triangle-free and 3-colorable.
        \item $\Xi$ contains a subset $X$ of four vertices such that in some proper 3-coloring of $\Xi$, the closed neighborhood $N(X) \setminus X$ is monochromatic, and $X$ satisfies one of the following:
        \begin{enumerate}
            \item $X$ contains at most two triangles, all triangles of $\Xi$ are contained within $X$, three vertices of $X$ have no neighbors outside of $X$, and $\Xi[X]$ is connected. Moreover, if a vertex in $X$ indeed has neighbors outside $X$, then it is adjacent to all other vertices in $X$.
            \item There exists a labeling of the vertices in $X$ by $\{a,b,c,d\}$, such that all triangles in $\Xi$ contain the $b-c$ edge, $a$ has no neighbors within $X$, $d$ has no neighbors outside of $X$, $N(c)\setminus\{b\}\subseteq N(b)$, and one of the following is satisfied:
            \begin{enumerate}[label = (\roman*)]
            \item $N(c)\setminus X \subseteq N(a)$
            \item $N(a) = \emptyset$.
            \end{enumerate}
        \end{enumerate}
    \end{enumerate}
\end{theorem}
\begin{proof}
    \textbf{Forward direction.} Let $G$ be an arbitrary $\PSL(2,13)$-solvable group. By \thref{PSLcoloring}, we know $\pgc(G)$ is 3-colorable. Thus, we may assume $G$ contains a triangle because otherwise $\pgc(G)$ matches Condition $1$. Similarly, by \thref{PSLsatisfies}, we may assume that $G$ satisfies the assumptions of \thref{structuralresult}. We also know by \thref{PSLedgenonexistence} that $3-p \not\in\pgc(G)$ for all $p\in \pi(G)\setminus\pi(\PSL(2,13))$. We break into cases according to the number of missing edges of $\pgc(\PSL(2,13))$. 
    
    \underline{Case 1: Zero edges of $\pgc(\PSL(2,13))$ are missing.} Since $2$ is not isolated from the primes of $\PSL(2,13)$, we know that $2$ and $7$ do not connect to primes outside of $\pi(\PSL(2,13))$ by \thref{2connectsoutward,PSL2noedge7noedge}. Therefore $2, 3$, and $7$ are isolated, only $13$—which is adjacent to $2$, $3$, and $7$—may connect outward, and $\pgc(G)$ contains exactly the two triangles within $\pgc(\PSL(2,13))$ by \thref{PSLtriangleedge}. Thus $\pgc(G)$ satisfies 2(a).
    
    \underline{Case 2: One edge of $\pgc(\PSL(2,13))$ is missing.} By \thref{PSLinternalstructure},$2-13 \in \pgc(G)$ and $3-13\in\pgc(G)$. Since $2$ is not isolated from the primes of $\PSL(2,13)$, we have by \thref{2connectsoutward,PSL2noedge7noedge} that $2$ and $7$ do not connect to primes outside of $\pi(\PSL(2,13))$. Therefore $2, 3$, and $7$ are isolated, only $13$—which is adjacent to $2$, $3$, and $7$—may connect outward, and the possible triangles are exactly one of $\{2,7,13\}$ or $\{3,7,13\}$ by \thref{PSLtriangleedge}. Thus $\pgc(G)$ satisfies 2(a).

    \underline{Case 3: Two edges edges of $\pgc(\PSL(2,13))$ are missing.} We consider subcases on the edges $2-7$ and $2-13$. 
    
    \textit{Case 3a: $2-7,2-13\in\pgc(G)$.} By \thref{2connectsoutward,PSL2noedge7noedge}, the primes $2$ and $7$ cannot connect to primes outside of $\pi(\PSL(2,13))$ and $3$ has no neighbors—since the edges $3-13$ and $3-7$ must have been the ones removed. So $\pgc(G)$ satisfies Condition 2(b)(ii) with the labeling $3\rightarrow a$, $7\rightarrow c$, $13\rightarrow b$, and $2\rightarrow d$. 
    
    \textit{Case 3b: Exactly one of $2-7$ or $2-13$ in $\pgc(G)$.} By \thref{2connectsoutward,PSL2noedge7noedge}, $2$ and $7$ cannot connect to primes outside of $\pi(\PSL(2,13))$. Then by \thref{PSLtriangleedge}, there cannot exist any triangles not contained within $\pi(\PSL(2,13))$. However, since an edge from each triangle within $\pi(\PSL(2,13))$ has been removed, this forms a contradiction with our assumption—so we do not need to classify this case. 
    
    \textit{Case 3c: $2-7, 2-13 \not\in\pgc(G)$.} We see that $N(2)\cap \pi(\PSL(2,13))=\emptyset$. First, suppose that $13$ is isolated from $\pi(G)\setminus \pi(\PSL(2,13))$. Then $\pgc(G)$ satisfies 2(b)(i) with the labeling $2\rightarrow a$, $7\rightarrow b$, $13 \rightarrow c$ and $3\rightarrow d$. Similarly, suppose $2$ is isolated from $\pi(G)\setminus\pi(\PSL(2,13))$ but $13$ is not. By \thref{PSL7subset}, $\pgc(G)$ satisfies 2(b)(ii) with the labeling $2\rightarrow a$, $13\rightarrow b$, $7 \rightarrow c$ and $3\rightarrow d$. Finally, suppose $2$ and $13$ both connect to $\pi(G)\setminus\pi(\PSL(2,13))$. By \thref{PSL7subset}, $\pgc(G)$ satisfies 2(b)(i) with the labeling $2\rightarrow a$, $13\rightarrow b$, $7 \rightarrow c$ and $3\rightarrow d$.
    
    \underline{Case 4: Three edges of $\pgc(\PSL(2,13))$ are missing.} Recall $7-13\in\pgc(G)$ by \thref{PSLtriangleedge}, so we proceed by considering possibilities on which other edge remains. If either $2-7\in\pgc(G)$ or $2-13\in\pgc(G)$ then by \thref{PSL2noedge7noedge}, we have that $7$ does not connect to $\pi(G)\setminus\pi(\PSL(2,13))$—which contradicts the existence of a triangle in $\pgc(G)$ by \thref{PSLtriangleedge}. The remaining edge also cannot be $3-7$ by \thref{PSLinternalstructure}. So we must only consider if the remaining edge is $3-13$. Then $N(2)\cap\pi(\PSL(2,13)) = \emptyset$ and 13 must connect to $\pi(G)\setminus\pi(\PSL(2,13))$ to form the triangle. By \thref{PSL7subset}, if $2$ does not connect to $\pi(G)\setminus\pi(\PSL(2,13))$ then $\pgc(G)$ satisfies 2(b)(ii) and if $2$ does connect to $\pi(G)\setminus\pi(\PSL(2,13))$ then $\pgc(G)$ satisfies 2(b)(i)—both with the labeling $2\rightarrow a$, $13\rightarrow b$, $7 \rightarrow c$ and $3\rightarrow d$.

    \underline{Case 5: Four edges of $\pgc(\PSL(2,13))$ are missing.}
    The only remaining edge must be $13-7$ by \thref{PSLtriangleedge}. So we have that $N(2)\cap\pi(\PSL(2,13)) = \emptyset$. Furthermore, 13 must connect to $\pi(G)\setminus\pi(\PSL(2,13))$ to form the triangle. By \thref{PSL7subset}, if $2$ does not connect to $\pi(G)\setminus\pi(\PSL(2,13))$, then $\pgc(G)$ satisfies 2(b)(ii), and if $2$ does connect to $\pi(G)\setminus\pi(\PSL(2,13))$, then $\pgc(G)$ satisfies 2(b)(i)—both with the labeling $2\rightarrow a$, $13\rightarrow b$, $7 \rightarrow c$ and $3\rightarrow d$. 
    
    These cases are exhaustive by \thref{PSLtriangleedge}, completing the forward direction.
    
    \textbf{Backward direction.} By \thref{solvableclassification} if $\Xi$ satisfies Condition $1$ then it is realizable by a solvable—hence $\PSL(2,13)$-solvable—group. Assume therefore that $\Xi$ contains a triangle. We proceed by cases once again.
    
    \underline{Case 1: $\Xi$ matches Condition 2(a).} By reviewing the graph structures in Appendix \hyperref[info]{A}, there exists a graph isomorphism $\varphi$ from $\Xi[X]$ to  $\pgc(E)$, where $E$ is either $\PSL(2,13)$ or $\Aut(\PSL(2,13))$. If no vertex in $X$ connects to $\Xi\setminus X$, the hypotheses of \cite[Lemma 2.3.6]{2023REU} are satisfied immediately by the trivial representation of $E$. Otherwise, let $v \in X$ connect outwards. By assumption $v$ is adjacent to all other vertices in $X$, so we may assume $\varphi(v)=13$. We now apply \cite[Lemma 2.3.6]{2023REU} with $E$ and its representations corresponding to the row vectors $[2,3,7]$ and $[2,3,7,13]$ in Appendix \hyperref[info]{A}.

    \underline{Case 2: $\Xi$ matches Condition 2(b)(i).} By assumption $\Xi$ contains a triangle, so there must be at least one edge in $\Xi[X]$. By reviewing the graph structures in Appendix \hyperref[info]{A}, there exists a graph isomorphism $\varphi$ from $\Xi[X]$ to $\pgc(E)$, where $E$ is one of $2.\PSL(2,13), \mathbb{F}_3^{36}\rtimes 2.\PSL(2,13)$, or $C_3\times 2.\PSL(2,13)$ such that $\varphi$ maps $a\rightarrow 2$, $b\rightarrow 13$, $c\rightarrow 7$, and $d\rightarrow 3$. We see by \thref{PSLedgenonexistence} and the assumption on the coloring of $\Xi$ that all cases for $E$ satisfy the hypotheses of \cite[Lemma 2.3.6]{2023REU} or \cite[Lemma 1.3]{2024PSL}—depending on if $E=2.\PSL(2,13)$ or $E\in\{\mathbb{F}_3^{36}\rtimes 2.\PSL(2,13), C_3\times 2.\PSL(2,13)\}$, respectively—if we show that $2.\PSL(2,13)$ has the correct fixed point information as presented in Appendix \hyperref[info]{A}. We consider all the allowed cases for $v \in \Xi\setminus X$. If $v\not\in N(a) \cup N(b) \cup N(c)$ then the representation for $2.\PSL(2,13)$ corresponding to $[2,3,7,13]$ suffices. If $v \in N(a)$ but not in $N(b)\cup N(c)$ then the representation for $[3,7,13]$ suffices. If $v\in N(b)$ but not in $N(a)\cup N(c)$ then the representation for $[2,3,7]$ suffices. If $v\in N(a)\cap N(b)$ but not in $N(c)$ then the representation for $[3,7]$ suffices. If $v\in N(c)$ then we also know $v\in N(b)\cap N(a)$ so the representation for $[3]$ suffices.

    \underline{Case 3: $\Xi$ matches Condition 2(b)(ii).} Generate $\Xi^\prime$ satisfying 2(b)(i) by adding the edge $a-v$ for all $v\in N(c)\setminus X$. Let $G$ be a $\PSL(2,13)$-solvable group realizing $\Xi^\prime$ according to case 2. Then $G\times C_2$ realizes $\Xi$.
\end{proof}
\begin{figure}[h!]
\centering
\begin{tikzpicture}[
  node/.style={circle, draw, minimum size=8mm},
  mid arrow/.style={
    decoration={markings, mark=at position 0.5 with {\arrow[scale=1.2]{>}}},
    postaction={decorate},
    line width=1.2pt,
    >=Stealth
  }
]

\node[node, fill =blue!30] (D) at (0,0) {d};
\node[node, fill =green!30] (C) at (2,0) {c};
\node[node, fill =red!30] (B) at (0,2) {b};
\node[node, fill =green!30] (A) at (2,2) {a};

\node[node, fill =blue!30] (E) at (5.5,4) {};
\node[node, fill =blue!30] (F) at (5.5,2) {};
\node[node, fill =blue!30] (G) at (5.5,0) {};
\node[node, fill =blue!30] (H) at (5.5,-2){};

\node[node, fill =red!30] (I) at (7.5,4) {};
\node[node, fill =red!30] (J) at (7.5,2) {};
\node[node, fill =red!30] (K) at (7.5,0) {};
\node[node, fill =red!30] (L) at (7.5,-2) {};

\node[node, fill =green!30] (M) at (9.5,4) {};
\node[node, fill =green!30] (N) at (9.5,2) {};
\node[node, fill =green!30] (O) at (9.5,0) {};
\node[node, fill =green!30] (P) at (9.5,-2) {};

\draw [line width=1.2pt](D) -- (B);
\draw [line width=1.2pt, color=orange](B) -- (C);
\draw [line width=1.2pt, color=orange](B)--(E);
\draw [line width=1.2pt] (A) -- (E);
\draw[line width=1.2pt, color=orange] (C) -- (E);
\draw[line width=1.2pt, color=orange] (C) -- (H);
\draw[line width=1.2pt, color=orange] (B) to[out =290, in = 180](H);
\draw[line width=1.2pt] (B)to[out=60, in = 120](F);
\draw[line width=1.2pt] (A) --(G);
\draw [line width=1.2pt] (A)--(H);

\draw[line width=1.2pt] (I)--(E);
\draw [line width=1.2pt] (I)--(F);
\draw [line width=1.2pt] (K)--(F);
\draw [line width=1.2pt] (L) -- (H);
\draw [line width=1.2pt] (M) -- (J);
\draw [line width=1.2pt] (I) -- (N);
\draw [line width=1.2pt] (J) -- (N);
\draw[line width=1.2pt] (J) -- (O);
\draw[line width=1.2pt] (F) -- (O);
\draw [line width=1.2pt] (P) -- (L);

\end{tikzpicture}
\caption{A graph satisfying Condition 2(b)(i) of \thref{pslclassification} with triangles highlighted.}
\end{figure}

\appendix
\renewcommand{\thesection}{Appendix \Alph{section} -\hspace{-12pt}}
\section{Group Information}
\label{info}
The basic information for the simple groups considered in this paper is recorded in the following table. The Schur multiplier, outer automorphisms, and prime divisors can all be found in \cite{ATLAS}. The prime graph complements are computed in \cite{GAP} by identifying element orders within the group.

\begin{center}
\begin{tabular}{|c|c|c|c|c|}
\hline
$G$ & Schur Multiplier & $\operatorname{Out}(G)$ & $\pi(G)$ & $\pgc(G)$\\
\hline

$A_{11}$ & $2$ & $2$ & $\{2,3,5,7,11\}$ &
$\vcenter{\hbox{\scalebox{0.75}{%
\begin{tikzpicture}[state/.style={circle, draw, minimum size=0.75cm, inner sep=0cm, outer sep=0cm}]
\node[state] (A) at (0,0) {$2$};
\node[state] (B) at (1.5,0) {$3$};
\node[state] (C) at (3,0) {$5$};
\node[state] (D) at (0.75,-1.5) {$7$};
\node[state] (E) at (2.25,-1.5) {$11$};

\draw[line width=1.2pt] (A)--(E);
\draw[line width=1.2pt] (B)--(E);
\draw[line width=1.2pt] (C)--(D);
\draw[line width=1.2pt] (C)--(E);
\draw[line width=1.2pt] (D)--(E);
\end{tikzpicture}}}}$ \\
\hline

$A_{12}$ & $2$ & $2$ & $\{2,3,5,7,11\}$ &
$\vcenter{\hbox{\scalebox{0.75}{%
\begin{tikzpicture}[state/.style={circle, draw, minimum size=0.75cm, inner sep=0cm, outer sep=0cm}]
\node[state] (A) at (0,0) {$2$};
\node[state] (B) at (1.5,0) {$3$};
\node[state] (C) at (3,0) {$5$};
\node[state] (D) at (0.75,-1.5) {$7$};
\node[state] (E) at (2.25,-1.5) {$11$};
\draw[line width=1.2pt] (A)--(E);
\draw[line width=1.2pt] (B)--(E);
\draw[line width=1.2pt] (C)--(E);
\draw[line width=1.2pt] (D)--(E);
\end{tikzpicture}}}}$ \\
\hline

$A_{13}$ & $2$ & $2$ & $\{2,3,5,7,11,13\}$ &
$\vcenter{\hbox{\scalebox{0.75}{%
\begin{tikzpicture}[state/.style={circle, draw, minimum size=0.75cm, inner sep=0cm, outer sep=0cm}]
\node[state] (A) at (.25,0) {$2$};
\node[state] (B) at (1.5,0) {$3$};
\node[state] (C) at (3,0) {$5$};
\node[state] (D) at (4.25,0) {$7$};
\node[state] (E) at (1.5,-1.5) {$11$};
\node[state] (F) at (3,-1.5) {$13$};
\draw[line width=1.2pt] (A)--(F);
\draw[line width=1.2pt] (A)--(E);
\draw[line width=1.2pt] (B)--(F);
\draw[line width=1.2pt] (B)--(E);
\draw[line width=1.2pt] (C)--(F);
\draw[line width=1.2pt] (C)--(E);
\draw[line width=1.2pt] (D)--(E);
\draw[line width=1.2pt] (D)--(F);
\draw[line width=1.2pt] (E)--(F);
\end{tikzpicture}}}}$ \\
\hline

$A_{14}$ & $2$ & $2$ & $\{2,3,5,7,11,13\}$ &
$\vcenter{\hbox{\scalebox{0.75}{%
\begin{tikzpicture}[state/.style={circle, draw, minimum size=0.75cm, inner sep=0cm, outer sep=0cm}]
\node[state] (A) at (.25,0) {$2$};
\node[state] (B) at (1.5,0) {$3$};
\node[state] (C) at (3,0) {$5$};
\node[state] (D) at (4.25,0) {$7$};
\node[state] (E) at (1.5,-1.5) {$11$};
\node[state] (F) at (3,-1.5) {$13$};
\draw[line width=1.2pt] (A)--(F);
\draw[line width=1.2pt] (A)--(E);
\draw[line width=1.2pt] (B)--(F);
\draw[line width=1.2pt] (C)--(F);
\draw[line width=1.2pt] (C)--(E);
\draw[line width=1.2pt] (D)--(E);
\draw[line width=1.2pt] (D)--(F);
\draw[line width=1.2pt] (E)--(F);
\end{tikzpicture}}}}$ \\
\hline

$A_{15}$ & $2$ & $2$ & $\{2,3,5,7,11,13\}$ &
$\vcenter{\hbox{\scalebox{0.75}{%
\begin{tikzpicture}[state/.style={circle, draw, minimum size=0.75cm, inner sep=0cm, outer sep=0cm}]
\node[state] (A) at (.25,0) {$2$};
\node[state] (B) at (1.5,0) {$3$};
\node[state] (C) at (3,0) {$5$};
\node[state] (D) at (4.25,0) {$7$};
\node[state] (E) at (1.5,-1.5) {$11$};
\node[state] (F) at (3,-1.5) {$13$};
\draw[line width=1.2pt] (A)--(F);
\draw[line width=1.2pt] (B)--(F);
\draw[line width=1.2pt] (C)--(F);
\draw[line width=1.2pt] (C)--(E);
\draw[line width=1.2pt] (D)--(E);
\draw[line width=1.2pt] (D)--(F);
\draw[line width=1.2pt] (E)--(F);
\end{tikzpicture}}}}$ \\
\hline

$A_{16}$ & $2$ & $2$ & $\{2,3,5,7,11,13\}$ &
$\vcenter{\hbox{\scalebox{0.75}{%
\begin{tikzpicture}[state/.style={circle, draw, minimum size=0.75cm, inner sep=0cm, outer sep=0cm}]
\node[state] (A) at (.25,0) {$2$};
\node[state] (B) at (1.5,0) {$3$};
\node[state] (C) at (3,0) {$5$};
\node[state] (D) at (4.25,0) {$7$};
\node[state] (E) at (1.5,-1.5) {$11$};
\node[state] (F) at (3,-1.5) {$13$};
\draw[line width=1.2pt] (A)--(F);
\draw[line width=1.2pt] (C)--(F);
\draw[line width=1.2pt] (D)--(E);
\draw[line width=1.2pt] (D)--(F);
\draw[line width=1.2pt] (E)--(F);
\end{tikzpicture}}}}$ \\
\hline

$\PSL(2,13)$ & $2$ & $2$ & $\{2,3,7,13\}$ &
$\vcenter{\hbox{\scalebox{0.75}{%
\begin{tikzpicture}[state/.style={circle, draw, minimum size=0.75cm, inner sep=0cm, outer sep=0cm}]
\node[state] (A) at (0,0) {$3$};
\node[state] (B) at (1.5,0) {$7$};
\node[state] (C) at (0,1.5) {$13$};
\node[state] (D) at (1.5,1.5) {$2$};
\draw [line width=1.2pt](A) -- (B);
\draw [line width=1.2pt](B) -- (D);
\draw [line width=1.2pt](A) -- (C);
\draw [line width=1.2pt](D) -- (C);
\draw [line width=1.2pt](B) -- (C);
\end{tikzpicture}}}}$ \\
\hline
$S_6(3)$ & $2$ & $2$ & $\{2,3,5,7,13\}$ &
$\vcenter{\hbox{\scalebox{0.75}{%
\begin{tikzpicture}[state/.style={circle, draw, minimum size=0.75cm, inner sep=0cm, outer sep=0cm}]
\node[state] (A) at (0,0) {$2$};
\node[state] (B) at (1.5,0) {$3$};
\node[state] (C) at (3,0) {$5$};
\node[state] (D) at (0.75,-1.5) {$7$};
\node[state] (E) at (2.25,-1.5) {$13$};
\draw[line width=1.2pt] (A)--(D);
\draw[line width=1.2pt] (B)--(E);
\draw[line width=1.2pt] (C)--(D);
\draw[line width=1.2pt] (C)--(E);
\draw[line width=1.2pt] (D)--(E);
\draw[line width=1.2pt] (B)--(D);
\end{tikzpicture}}}}$ \\
\hline$HS$ & $2$ & $2$ & $\{2,3,5,7,11\}$ &
$\vcenter{\hbox{\scalebox{0.75}{%
\begin{tikzpicture}[state/.style={circle, draw, minimum size=0.75cm, inner sep=0cm, outer sep=0cm}]
\node[state] (A) at (0,0) {$2$};
\node[state] (B) at (1.5,0) {$3$};
\node[state] (C) at (3,0) {$5$};
\node[state] (D) at (0.75,-1.5) {$7$};
\node[state] (E) at (2.25,-1.5) {$11$};
\draw[line width=1.2pt] (A)--(D);
\draw[line width=1.2pt] (A)--(E);
\draw[line width=1.2pt] (B)--(E);
\draw[line width=1.2pt] (C)--(D);
\draw[line width=1.2pt] (C)--(E);
\draw[line width=1.2pt] (D)--(E);
\draw[line width=1.2pt] (B)--(D);
\end{tikzpicture}}}}$ \\
\hline
$G_2(4)$ & $2$ & $2$ & $\{2,3,5,7,13\}$ &
$\vcenter{\hbox{\scalebox{0.75}{%
\begin{tikzpicture}[state/.style={circle, draw, minimum size=0.75cm, inner sep=0cm, outer sep=0cm}]
\node[state] (A) at (0,0) {$2$};
\node[state] (B) at (1.5,0) {$3$};
\node[state] (C) at (3,0) {$5$};
\node[state] (D) at (0.75,-1.5) {$7$};
\node[state] (E) at (2.25,-1.5) {$13$};
\draw[line width=1.2pt] (A)--(D);
\draw[line width=1.2pt] (A)--(E);
\draw[line width=1.2pt] (B)--(E);
\draw[line width=1.2pt] (C)--(D);
\draw[line width=1.2pt] (C)--(E);
\draw[line width=1.2pt] (D)--(E);
\end{tikzpicture}}}}$ \\
\hline

\end{tabular}
\end{center}

The following groups have the given fixed point information, calculated as in \cite[Table 3]{2023REU} with \cite{GAP,CTblLib}. A row vector $[p_1, p_2, \ldots, p_n]$ appears in this table if and only if there exists a complex irreducible representation of $G$ such that, for each $1\leq i\leq n$, some $x\in G$ of order $p_i$ acts with a fixed point and, for every other prime $q\in\pi(G)$, if $y\in G$ has order $q$ then it acts Frobeniusly.
\begin{center}
    \begin{tabular}[h]{|c|c|}
    \hline
        $G$&  Fixed Point Information\\
        \hline
        $A_{11}$ & $ [2,3,5,7], [2,3,5,7,11] $ \\
        $2.A_{11}$ & $ [2,3,5,7],[2,3,5,7,11] $ \\
        $A_{12}$ & $ [2,3,5,7,11] $ \\
        $2.A_{12}$ & $ [2,3,5,7,11] $ \\
        $A_{13}$ & $ [2,3,5,7,11], [2,3,5,7,11,13] $ \\
        $2.A_{13}$ & $ [2,3,5,7,11], [2,3,5,7,11,13] $ \\
        $\PSL(2,13)$ & $ [ 2, 3, 7 ], [ 2, 3, 7, 13 ] $ \\
        $\Aut(\PSL(2,13))$ & $ [ 2, 3, 7 ], [ 2, 3, 7, 13 ] $ \\
        $2.\PSL(2,13)$ & $ [ 2, 3, 7 ], [ 2, 3, 7, 13 ], [ 3 ], [ 3, 7 ], [ 3, 7, 13 ] $ \\
        $S_6(3)$&$[ 2, 3, 5, 7, 13 ]$\\
        $2.S_6(3)$& $[ 2, 3, 5, 7, 13 ]$\\
        $HS$&$[ 2, 3, 5, 7, 11 ]$\\
        $2.HS$& $[ 2, 3, 5, 7, 11 ]$\\
        $G_2(4)$&$[ 2, 3, 5, 7, 13 ]$\\
        $2.G_2(4)$& $[ 2, 3, 5 ], [ 2, 3, 5, 7, 13 ]$\\
        \hline
    \end{tabular}
\end{center}

Finally, we record some additional information about some $\PSL(2,13)$-solvable groups which is needed in Subsection \ref{PSL} to prove \thref{pslclassification}. 

Most of the following prime graph complements can be directly computed in \cite{GAP} by identifying element orders within the group. The prime graph complement for the semidirect product group is instead found in \cite{GAP} by computing some irreducible representations for $2.\PSL(2,13)$ over $\mathbb{F}_3$ and identifying a representation where all elements of order $13$ act Frobeniusly but there exist elements of order $7$ which do not.
\begin{center}
\begin{tabular}{|c|c|}
\hline
$G$ & $\pgc(G)$\\
\hline

$2.\PSL(2,13)$ &
$\vcenter{\hbox{\scalebox{0.75}{%
\begin{tikzpicture}[state/.style={circle, draw, minimum size=0.75cm, inner sep=0cm, outer sep=0cm}]
\node[state] (A) at (0,0) {$3$};
\node[state] (B) at (1.5,0) {$7$};
\node[state] (C) at (0,1.5) {$13$};
\node[state] (D) at (1.5,1.5) {$2$};
\draw[line width=1.2pt] (A) -- (B);
\draw[line width=1.2pt] (A) -- (C);
\draw[line width=1.2pt] (B) -- (C);
\end{tikzpicture}}}}$ \\
\hline

$\Aut(\PSL(2,13))$ &
$\vcenter{\hbox{\scalebox{0.75}{%
\begin{tikzpicture}[state/.style={circle, draw, minimum size=0.75cm, inner sep=0cm, outer sep=0cm}]
\node[state] (A) at (0,0) {$3$};
\node[state] (B) at (1.5,0) {$7$};
\node[state] (C) at (0,1.5) {$13$};
\node[state] (D) at (1.5,1.5) {$2$};
\draw[line width=1.2pt] (A) -- (B);
\draw[line width=1.2pt] (A) -- (C);
\draw[line width=1.2pt] (D) -- (C);
\draw[line width=1.2pt] (B) -- (C);
\end{tikzpicture}}}}$ \\
\hline

$C_3 \times 2.\PSL(2,13)$ &
$\vcenter{\hbox{\scalebox{0.75}{%
\begin{tikzpicture}[state/.style={circle, draw, minimum size=0.75cm, inner sep=0cm, outer sep=0cm}]
\node[state] (A) at (0,0) {$3$};
\node[state] (B) at (1.5,0) {$7$};
\node[state] (C) at (0,1.5) {$13$};
\node[state] (D) at (1.5,1.5) {$2$};
\draw[line width=1.2pt] (B) -- (C);
\end{tikzpicture}}}}$ \\
\hline

$\mathbb{F}_3^{36}\rtimes 2.\PSL(2,13)$ &
$\vcenter{\hbox{\scalebox{0.75}{%
\begin{tikzpicture}[state/.style={circle, draw, minimum size=0.75cm, inner sep=0cm, outer sep=0cm}]
\node[state] (A) at (0,0) {$3$};
\node[state] (B) at (1.5,0) {$7$};
\node[state] (C) at (0,1.5) {$13$};
\node[state] (D) at (1.5,1.5) {$2$};
\draw[line width=1.2pt] (A) -- (C);
\draw[line width=1.2pt] (B) -- (C);
\end{tikzpicture}}}}$ \\
\hline

\end{tabular}
\end{center}

The information in the following table, calculated with \cite{GAP, CTblLib}, is to facilitate easy applications of \thref{suzgen}. A row vector $[p_1, p_2, \ldots, p_n]$ appears in the table if and only if there exists some $\chi\in\operatorname{IBr}_p(G)$ such that $\{p_1, p_2, \ldots, p_n\} = \{q\in\pi(G)\setminus\{p\}\mid\exists g \in G \operatorname{where} |g| = q, \frac{1}{|g|}\sum_{x\in\langle g \rangle}\chi(g)>0\}$.
\begin{center}
    \begin{tabular}{|c|c|c|}
        \hline
            $G$ & $p$ & Brauer Information \\
            \hline
            \multirow{4}{*}{$\PSL(2,13)$}& 2 & $[ 3 ], [ 3, 7 ], [ 3, 7, 13 ]$\\
            
            & 3 &$[ 2, 7 ], [ 2, 7, 13 ]$\\
            
            & 7 & $[ 2, 3 ], [ 2, 3, 13 ]$\\
            
            & 13 & $[ 2, 3, 7 ]$\\
           \hline
            \multirow{4}{*}{$2.\PSL(2,13)$}& 2 & $[ 3 ], [ 3, 7 ], [ 3, 7, 13 ]$\\
            
            & 3 &$[ \ ], [ 2, 7 ], [ 2, 7, 13 ], [ 7 ], [ 7, 13 ]$\\
            
            & 7 & $[ 2, 3 ], [ 2, 3, 13 ], [ 3 ], [ 3, 13 ]$\\
            
            & 13 & $[ \ ], [ 2, 3, 7 ], [ 3 ], [ 3, 7 ]$\\
            \hline
    \end{tabular}
\end{center}
\section{GAP Computations}
The computations performed in GAP \cite{GAP} and the GAP Character Table Library \cite{CTblLib} were fairly elementary, but for any interested reader we have included on our \href{https://github.com/lucasalland/Code-For-Prime-Graph-Complements-Calculations}{Github} the code used for computing all the information found in Appendix \hyperref[info]{A}.

\vspace{1cm}
\begin{center}
{\bf Acknowledgements}\\

\end{center}
\vspace{3pt}
This research was conducted under the NSF-REU grants DMS-2150205 and DMS-2447229 under the 
mentorship of the third author. The authors would like to thank Alexa Renner and Lixin Zheng
from the 2024 and 2023 Texas State REU Teams (respectively) for their help. They also gratefully acknowledge the financial support of NSF and thank Texas
State University for providing a great working environment. The third author was also partially 
supported by a grant from the Simons Foundation (SFI-MPS-TSM-00014153 to Thomas M. Keller).


\begin{thebibliography}{10}

\bibitem{CTblLib}
Thomas Breuer.
\newblock The \textsf{GAP} {C}haracter {T}able {L}ibrary, {V}ersion 1.3.1.
\newblock http://www.math.rwth-aachen.de/\~{}Thomas.Breuer/ctbllib, April 2020.
\newblock \textsf{GAP} package.

\bibitem{maslova}
Peter~J. Cameron and Natalia~V. Maslova.
\newblock Criterion of unrecognizability of a finite group by its
  {G}ruenberg-{K}egel graph.
\newblock {\em J. Algebra}, 607:186--213, 2022.

\bibitem{frattini}
Klaus Doerk and Trevor O. Hawkes.
\newblock {\em Finite Soluble Groups}.
\newblock De Gruyter Expositions in Mathematics. W. de Gruyter, 1992.

\bibitem{2022REU}
Timothy~J. Edwards, Thomas Michael Keller, Ryan~M. Pesak, and
  Karthik~Sellakumaran Latha.
\newblock The prime graphs of groups with arithmetically small composition
  factors.
\newblock {\em Ann. Mat. Pura Appl. (4)}, 203(2):945--973, 2024.

\bibitem{Flavell}
Paul Flavell.
\newblock A {H}all-{H}igman-{S}hult type theorem for arbitrary finite groups.
\newblock {\em Invent. Math.}, 164(2):361--397, 2006.

\bibitem{GAP}
The GAP~Group.
\newblock {\em {GAP -- Groups, Algorithms, and Programming, Version 4.14.0}},
  2024.

\bibitem{2015REU}
Alexander Gruber, Thomas~Michael Keller, Mark~L. Lewis, Keeley Naughton, and
  Benjamin Strasser.
\newblock A characterization of the prime graphs of solvable groups.
\newblock {\em J. Algebra}, 442:397--422, 2015.

\bibitem{GraphTheory}
Lih-Hsing Hsu and Cheng-Kuan Lin.
\newblock {\em Graph Theory and Interconnection Networks}.
\newblock CRC Press, 2008.

\bibitem{2021REU}
Ziyu Huang, Thomas~Michael Keller, Shane Kissinger, Wen Plotnick, Maya Roma,
  and Yong Yang.
\newblock A classification of the prime graphs of pseudo-solvable groups.
\newblock {\em J. Group Theory}, 27(1):89--117, 2024.

\bibitem{HuppertEnglish}
Bertram Huppert.
\newblock {\em Finite Groups I}.
\newblock Springer International Publishing, 2025.

\bibitem{IsaacsFiniteGroups}
{I. Martin} Isaacs.
\newblock {\em Finite Group Theory}.
\newblock Graduate studies in mathematics. American Mathematical Society, 2008.

\bibitem{Suz}
{Thomas Michael} Keller, Zachary Martin, Alexa Renner, Gabriel Roca, and Eric
  Yu.
\newblock Classification of the prime graphs of $\operatorname{Sz}(8)$-,
  $\operatorname{Sz}(32)$-, and $\operatorname{PSL}(2, 2^5)$-solvable groups, submitted, 2024. 	arXiv:2410.21063 [math.GR]
\bibitem{2024PSL}
Thomas~Michael Keller, Zachary Martin, Alexa Renner, Gabriel Roca, and Eric Yu.
\newblock Criteria for classifying prime graphs of {PSL}(2, $q$)-solvable groups, submitted, 2025. arXiv:2510.21979 [math.GR]

\bibitem{2023REU}
Thomas~Michael Keller, Gavin Pettigrew, Saskia Solotko, and Lixin Zheng.
\newblock Classifying prime graphs of finite groups – a methodical approach.
\newblock {\em J. Pure Appl. Algebra}, 229(11):108089, 40pp., 2025.

\bibitem{Serre}
Jean-Pierre Serre.
\newblock {\em Linear Representations of Finite Groups}.
\newblock Graduate Texts in Mathematics. Springer, 1st edition, 1996.

\bibitem{aurelstackexchange}
Aurel (\url{https://mathoverflow.net/users/40821/aurel}).
\newblock A question on groups having a subgroup which fixes a vector in every
  irreducible representations.
\newblock MathOverflow.
\newblock \url{https://mathoverflow.net/q/479910} (version: 2024-10-02).

\bibitem{jashastackexchange}
Jasha (\url{https://math.stackexchange.com/users/163845/jasha}).
\newblock Proof that a group representation matrix is diagonalizable?
\newblock Mathematics Stack Exchange.
\newblock \url{https://math.stackexchange.com/q/969624} (version: 2014-10-12).

\bibitem{Webb}
Peter Webb.
\newblock {\em A Course in Finite Group Representation Theory}.
\newblock Cambridge studies in advanced mathematics. Cambridge University Press,
  2016.

\bibitem{Williams}
J.~S. Williams.
\newblock Prime graph components of finite groups.
\newblock {\em J. Algebra}, 69(2):487--513, 1981.

\bibitem{ATLAS}
Robert~A. Wilson, Peter Walsh, Jennifer Tripp, Ibrahim Suleiman, Richard
  Parker, Simon Norton, Steve Nickerson, Stephen Linton, John Bray, and Richard
  Abbott.
\newblock Atlas of finite group representations — version 3.
\newblock \url{https://brauer.maths.qmul.ac.uk/Atlas/v3/}.

\end{thebibliography}
\end{document}